\documentclass[11pt, letterpaper, leqno]{amsart} 

\usepackage{graphicx} 
\usepackage{amsmath,amsthm,amssymb} 
\usepackage{amsfonts} 
\usepackage{dsfont}
\usepackage{mathrsfs} 
\usepackage[T1]{fontenc} 
\usepackage[utf8]{inputenc} 

\usepackage[dvipsnames, svgnames]{xcolor} 
\usepackage[shortlabels]{enumitem} 
\usepackage[colorlinks,citecolor=citec,hypertexnames=false,urlcolor=urlc,breaklinks=true, pagebackref]{hyperref} 

\usepackage{accents} 

\usepackage{amsbsy} 
\usepackage{amscd} 
\usepackage{latexsym} 
\usepackage{txfonts} 

\usepackage{exscale} 
\usepackage{bbm} 
\usepackage{mathtools} 
\usepackage{bbding}

\usepackage{geometry} 

\makeatletter
\providecommand\@dotsep{4.5}
\makeatother

\usepackage{thmtools} 

\newif\ifsoldark
\newif\ifsollight
\newif\ifclassic
\newif\ifplain

\usepackage{cite} 

\renewcommand*\backref[1]{\ifx#1\relax \else (p. #1) \fi} 

\numberwithin{equation}{section} 

\theoremstyle{plain}
\newtheorem{theorem}[equation]{Theorem}
\newtheorem{lemma}[equation]{Lemma}
\newtheorem{corollary}[equation]{Corollary}
\newtheorem{proposition}[equation]{Proposition}

\theoremstyle{definition}
\newtheorem{definition}[equation]{Definition}

\newtheorem{problem}[equation]{Problem}
\newtheorem{conjecture}[equation]{Conjecture}
\newtheorem{example}[equation]{Example}
\newtheorem{assumption}[equation]{Assumption}

\theoremstyle{remark}
\newtheorem{remark}[equation]{Remark}

\newcommand{\dist}{\operatorname{dist}}

\newcommand{\rad}{\operatorname{rad}}

\newcommand{\dv}{\operatorname{div}}

\newcommand{\n}[1]{\mathscr{#1}}
\newcommand{\bb}[1]{\mathbb{#1}}

\newcommand{\RNum}[1]{\uppercase\expandafter{\romannumeral #1\relax}}

\newcommand{\Tr}{\operatorname{Tr}}

\DeclareMathOperator{\supp}{supp}

\def\Xint#1{\mathchoice
	{\XXint\displaystyle\textstyle{#1}}%
	{\XXint\textstyle\scriptstyle{#1}}%
	{\XXint\scriptstyle\scriptscriptstyle{#1}}%
	{\XXint\scriptscriptstyle\scriptscriptstyle{#1}}%
	\!\int}
\def\XXint#1#2#3{{\setbox0=\hbox{$#1{#2#3}{\int}$}
		\vcenter{\hbox{$#2#3$}}\kern-.5\wd0}}

\def\dashint{\Xint-}

\def\YYint#1#2#3{{\setbox0=\hbox{$#1{#2#3}{\iint}$}
		\vcenter{\hbox{$#2#3$}}\kern-.51\wd0}}
 
\newcommand{\ep}{\varepsilon}
\newcommand{\ra}{\rightarrow}
\newcommand{\lra}{\longrightarrow}
\newcommand{\m}[1]{\mathcal{#1}}

\newcommand{\f}[1]{\mathfrak{#1}}
\newcommand{\vertiii}[1]{{\left\vert\kern-0.25ex\left\vert\kern-0.25ex\left\vert #1 
		\right\vert\kern-0.25ex\right\vert\kern-0.25ex\right\vert}}

\newcommand{\loc}{\operatorname{loc}}

\newcommand*{\dt}[1]{%
	\accentset{\mbox{\Large\bfseries .}}{#1}}

\allowdisplaybreaks

\setlist{nosep} 

\colorlet{citec}{blue}
\colorlet{urlc}{blue}
\colorlet{toc}{blue}
\colorlet{hyperc}{blue}
\colorlet{bpcolor}{NavyBlue}
\colorlet{smcolor}{purple}
\colorlet{impcolor}{blue}
\colorlet{eqcolor}{black}
\colorlet{lmcolor}{black}
\colorlet{propcolor}{black}
\colorlet{thmcolor}{black}
\colorlet{defcolor}{black}
\colorlet{rmcolor}{black}
\colorlet{excolor}{black}

\usepackage{everypage} 
\usepackage{tikzpagenodes} 

\newcommand{\1}{{\mathds 1}}

 \renewcommand{\a}{{\bf a}}
\newcommand{\B}{{\bf B}}

\begin{document}

\author[B. Poggi]{Bruno Poggi}
\email[]{poggi@mat.uab.cat}
\address{Departament de  Matemàtiques, Universitat Autònoma de Barcelona, Bellaterra,  Catalunya}

\title[Applications of the landscape function]{Applications of the landscape function     for Schr\"odinger operators with singular potentials and irregular magnetic fields}
\date{\today}

\thanks{The results in this article were obtained while the author was partly supported by the University of Minnesota Doctoral Dissertation Fellowship,    the European Research Council (ERC) under the European Union's Horizon 2020 research and innovation programme (grant agreement 101018680), the Simons Collaborations in MPS 563916, SM, the NSF DMS grant 1344235, and the NSF DMS grant DMS-1839077. The author would like to thank Professor S. Mayboroda for many insightful and inspiring conversations.}

\begin{abstract} We resolve both a conjecture and  a problem of Z. Shen from the  90's regarding non-asymptotic   bounds on the eigenvalue counting function  of the magnetic Schr\"odinger operator  $L_{\a,V}=-(\nabla-i{\bf a})^2+V$ with a singular or irregular magnetic field $\B$ on $\bb R^n$, $n\geq3$. We do this by constructing a new landscape function for $L_{\a,V}$, and proving its corresponding   uncertainty principle,  under  certain directionality assumptions on $\B$, but with no assumption on $\nabla\B$.   These results arise as applications of our study of the Filoche-Mayboroda landscape function $u$, a solution to  the equation $L_Vu=-\dv A\nabla u+Vu=1$, on unbounded Lipschitz domains in $\bb R^n$, $n\geq1$, and   $0\leq V\in L^1_{\loc}$, under a mild decay   condition on the Green's function. For  $L_V$,  we prove a priori exponential decay of Green's function,  eigenfunctions, and Lax-Milgram solutions in an Agmon distance with weight $1/u$, which may degenerate. Similar a priori results hold for $L_{\a,V}$. Furthermore, when $n\geq3$ and $V$ satisfies a scale-invariant Kato condition and a weak doubling property, we show that $1/\sqrt u$ is pointwise equivalent to the   Fefferman-Phong-Shen maximal function $m(\cdot,V)$; in particular this gives a strong scale-invariant Harnack inequality for $u$, and a setting where the Agmon distance with weight $1/u$ is not too degenerate.  Finally, we  extend results from the literature for $L_{\a,V}$ regarding exponential decay of the fundamental solution and eigenfunctions,   to the  situation of   irregular magnetic fields with directionality assumptions.
\end{abstract}

\maketitle

{
	\hypersetup{linkcolor=toc}
	\tableofcontents
}
\hypersetup{linkcolor=hyperc}

\section{Introduction}

In the present paper, we investigate the link between the Filoche-Mayboroda landscape function for the Dirichlet problem and the exponential decay of solutions to Schr\"odinger operators 
\begin{equation}\label{eq.s}
	L=-\dv A\nabla+V
\end{equation}
on a (possibly unbounded) domain $\Omega\subseteq\bb R^n$, $n\geq1$, with non-negative potentials $V\in L^1_{\loc}(\Omega)$ and $A=(A_{ij})$ a (not necessarily symmetric) matrix of bounded measurable coefficients on $\Omega$ verifying the uniform ellipticity condition
\begin{equation}\label{eq.elliptic}
	\lambda|\xi|^2\leq\sum_{i,j=1}^nA_{ij}(x)\xi_i\xi_j,\qquad \Vert A\Vert_{L^{\infty}(\Omega)}\leq\frac1{\lambda},\quad x\in\Omega,\quad\xi\in\bb R^n,
\end{equation}  
for some $\lambda\in(0,1)$. Here, the landscape function is the solution $u$ to the problem $Lu=1$ with Dirichlet boundary conditions. A part of this story (with Neumann boundary conditions) was already considered for eigenfunctions in \cite{adfjm}, and partly our goal here will be to further study the exponential decay as it pertains to Lax-Milgram solutions, Dirichlet eigenfunctions, and Green's functions (Theorems \ref{thm.lmr}, \ref{thm.eigen}, \ref{thm.greenexp}), while also removing the boundedness assumption on $V$ and on the domain. We will show that, under a mild decay condition  on the Green's function, a landscape function exists over unbounded domains (Theorem \ref{thm.exist}), and that it can be used to study a priori exponential decay of the Green's function. Moreover, we will show that, under natural assumptions, an appropriate power of this landscape function is pointwise equivalent to the Fefferman-Phong-Shen maximal function (\ref{eq.fpsm}) (Theorem \ref{thm.compare}), which is a tool that has been used to study many related properties of the Schr\"odinger operators, indicating that the landscape function is in some sense both a generalization and a regularization of this maximal function. In particular, in this setting $u$ verifies certain interesting properties (Corollary \ref{cor.u}) such as a strong scale-invariant Harnack inequality. Finally, we will apply these results to   exhibit the existence of a new  landscape function (Theorem \ref{thm.magup}) related to the magnetic Schr\"odinger operator 
\begin{equation}\label{eq.ms}
	L_{{\bf a},V}=-(\nabla-i{\bf a})^2+V,
\end{equation}
and we will use this landscape function, with its  uncertainty principle, to also obtain a priori exponential decay for this operator (Theorem \ref{thm.expdecayms}), and recover many results (Corollary \ref{cor.upm} and Section \ref{sec.cor})  for this operator under  directionality assumptions (Section \ref{sec.main2}) on the \emph{magnetic field}
\begin{equation}\label{eq.magfield}
	{\bf B}(x):=\operatorname{curl}{\bf a}(x)=(b_{jk}(x))_{1\leq j,k\leq n}=\Big(\frac{\partial a_j(x)}{\partial x_k}-\frac{\partial a_k(x)}{\partial x_j}\Big)_{1\leq j, k\leq n},
\end{equation}
which in particular allow for non-positive singular potentials $V$ and completely relax usual conditions on $\nabla\B$ found in related literature, solving a conjecture (Conjecture \ref{conj.s1}) and an open problem (Problem \ref{pro.s2}) of Zhongwei Shen from the 90's in the process.

Before precisely stating our main results in Sections \ref{sec.main} and \ref{sec.main2a}, let us provide the appropriate historical context and define the meaningful notions.

In 2012, a simple but remarkably versatile tool was introduced by Filoche and Mayboroda \cite{fm, fm2} to study the localization of eigenfunctions of Schr\"odinger operators. For the Schr\"odinger operator $L=-\Delta +V$, they considered a solution $u$ to the equation $Lu=1$, and showed that this function predicts the shape and location of localized low energy eigenfunctions of $L$ whose localization is caused either by the disorder of the potential or the geometry of the underlying domain. This solution has come to be known as the \emph{landscape function}, and over the last decade it has seen numerous appearances in mathematics \cite{st, lst, adfjm, adfjm2, dfm4, wz21, cwz, afmwz}, and theoretical and experimental physics \cite{fmphys, fmphys2, fmphys3, fmphys4, dmzdawf}, mainly concerning the study of the spectrum and eigenfunctions of Schr\"odinger operators acting under a disordered potential. In particular, it has yielded striking new results in the study of Anderson localization \cite{adfjm}, and a substantial improvement over the Weyl law \cite{dfm4}, for Schr\"odinger operators, although the numerical evidence \cite{adfjm2} suggests that the proven results so far capture only a fraction of its predictive power. The concept of the landscape function has been used effectively to study Anderson-type localization for Schr\"odinger operators on quantum graphs \cite{hm2020}, for the tight binding model \cite{afmwz}, for the case of many-body localization \cite{blg2020}, and even for random $M-$matrices \cite{fmt}. 

Finding a landscape function for the magnetic Schr\"odinger operator has remained an open question, and is of independent interest. Main difficulties   have been  that solutions to $L_{\a,V}u=f$ are generally complex-valued, and that the algebra of the equation has not been favorable towards finding an analogous conjugation (\ref{eq.conj-intro}) that would yield a useful \emph{uncertainty principle} (\ref{eq.up0}).  In Section \ref{sec.main2a}, we exhibit a new candidate for a landscape function for $L_{\a,V}$, which as we will see, has many  a priori properties similar to those of the landscape function of Filoche-Mayboroda.

The power of the landscape function $u$ to give so many mathematical and numerical results concerning Anderson localization (that is, localization of eigenfunctions via potentials that exhibit randomness) has been   referred to as ``almost magical'' \cite{fmt}, because even though the landscape function is \emph{deterministic}, it nonetheless can be used to study the effects of randomness on eigenfunctions and the integrated density of states. Moreover, the mathematical results \cite{adfjm} regarding the exponential decay of eigenfunctions for the Schr\"odinger operators with random potentials are phrased in terms of an Agmon-type decay with weight $\frac1u$, but the relationship of this Agmon distance  with the Euclidean distance has not been studied, and there may be no link when $V$ is too degenerate (that is, too close to $0$ too often). For instance, if the potential $V$ is $0$ outside of a small ball, then the Agmon distance with weight $1/u$ (with Neumann boundary conditions) between any two points $x,y\in\Omega$ is bounded above by a uniform constant regardless of how large $|x-y|$ is. Nevertheless, the numerical results predict   localization in the Euclidean sense when the potential is random. Since the landscape function is a deterministic tool, we might expect that, at least in some classical cases where the potential $V$ is not random, not degenerate, and has some regularity,   we should be able to connect the Agmon distance with weight $1/u$ to the Euclidean distance.

A classical  setting which has been studied extensively for its connection to harmonic analysis is the case of non-negative potentials $V$ in dimension $n\geq3$ which satisfy the \emph{reverse H\"older} inequality
\begin{equation}\label{eq.rh}
\Big(\dashint_BV^{\frac n2}\Big)^{\frac2n}\lesssim\dashint_B V,\qquad\text{for all balls }B.
\end{equation}
Such functions are said to belong to the space $RH_{n/2}$, and the best constant in (\ref{eq.rh}) is the \emph{$RH_{n/2}$ characteristic} of $V$. Examples of potentials in $RH_{n/2}$ include, for instance, all non-negative polynomials, and the functions $|x|^{\alpha}$ with $\alpha>-2$, while non-examples are potentials given by random variables (this is an easy exercise using the Borel-Cantelli lemma) and also exponentially decaying potentials. In the setting when (\ref{eq.rh}) holds, Z. Shen \cite{shenp} proved that the operators $\nabla(-\Delta+V)^{-1/2}$, $(-\Delta+V)^{-1/2}\nabla$, and $\nabla(-\Delta+V)^{-1}\nabla$ are bounded in $L^p$, extending similar $L^p$ estimates for non-negative polynomials of \cite{zhong}. A key idea  to prove these results is to first show that the \emph{fundamental solution} $\Gamma$ of the operator $-\Delta+V$ has polynomial decay of sufficiently large degree. In order to do this, Shen \cite{shenn} defined the maximal function $m(x,w)$, given by
\begin{equation}\label{eq.fpsm}
\frac1{m(x,w)}:=\sup\Big\{r\,:\,\frac1{r^{n-2}}\int_{B(x,r)}w(y)\,dy\leq C_1\Big\},
\end{equation}
which generalizes a smooth weight function  suited to work with non-negative polynomials \cite{smith91, hm, hn2, zhong},
and proved the following variant of the Fefferman-Phong  uncertainty principle  \cite{fefferman}
\begin{equation}\label{eq.fpu}
\int_{\bb R^n}m^2(x,V)f^2(x)\,dx\lesssim\int_{\bb R^n}\Big[|\nabla f(x)|^2+V(x)f^2(x)\Big]\,dx,\qquad\text{for all } f\in C_c^1(\bb R^n).
\end{equation}
The estimate (\ref{eq.fpu}) was used in \cite{shenn, shenp} to show the required polynomial decay of the fundamental solution. On the other hand, when $n=2$, M. Christ \cite{christ91} showed that the fundamental solution for the operator $-\Delta+V$ under certain related assumptions verifies quite sharp exponential decay estimates. The sharp exponential decay estimates for the fundamental solution in $n\geq3$ were later obtained in \cite{shenf} in terms of the geodesic $\rho(x,y,m^2)$ on the modified Agmon metric
\begin{equation}\label{eq.agmonmetric}
	ds^2=m(x,V)\big\{dx_1^2+dx_2^2+\ldots+dx_n^2\big\}.
\end{equation}
Thus, in the case when $n\geq3$ and $V\geq0$ verifies certain general assumptions (see below) which   imply (\ref{eq.rh}), Z. Shen proved the existence of positive constants $C$, $\ep_1$, and $\ep_2$, such that
\begin{equation}\label{eq.sharpexp}
\frac1{C|x-y|^{n-2}}e^{-\ep_1\rho(x,y,m^2)}\leq\Gamma(x,y)\leq\frac{C}{|x-y|^{n-2}}e^{-\ep_2\rho(x,y,m^2)},\qquad\text{for all }x,y\in\bb R^n, x\neq y,
\end{equation}
with the help of the uncertainty principle (\ref{eq.fpu}) and an Agmon-type argument \cite{agmon}. This sharp exponential decay result would later be extended by S. Mayboroda and the author of this article to encompass the case of Schr\"odinger operators $-\dv A\nabla+V$ with $A$ an elliptic, not necessarily symmetric matrix of bounded measurable coefficients \cite{mp}. 

In addition to the exponential decay of fundamental solutions \cite{shenf, mp} and the $L^p$ estimates for Schr\"odinger operators \cite{shenp, ab}, the maximal function $m$ (or  very similar ideas stemming from the Fefferman-Phong inequality) has also been used to study the Neumann problem for the Schr\"odinger operator $-\Delta+V$ \cite{shenn, mt}, bounds on the eigenvalue counting function \cite{shene, sheno, shenb}, heat kernels \cite{kurata}, and several analogous results concerning the magnetic Schr\"odinger operator $L_{\a,V}$ for $n\geq3$ \cite{shene, shenm, sheno, ks, kurata, kurata2000, benali, benali2, mp}. For these   results on $L_{\a,V}$, the authors  consider the maximal function $m(x,|\B|+V)$ where
\begin{equation}\label{eq.magb}
|\B|(x):=\sum_{j<k}|b_{jk}(x)|,\qquad x\in\bb R^n,
\end{equation}
and they naturally assume that $|\B|+V$ verifies (\ref{eq.rh}), to ensure that $m(\cdot,|\B|+V)$ has nice properties. However, in order to prove an uncertainty principle for the magnetic case, the following rather strong assumption is made on the gradient of the magnetic field
\begin{equation}\label{eq.gradb}
|\nabla\B(x)|\lesssim m(x,|\B|+V)^3,
\end{equation}
on all of the aforementioned works\footnote{We point out that in \cite{benali2}, the author does not assume (\ref{eq.gradb}), but still assumes a pointwise estimate on the gradient of the magnetic field.}. When $V$ is not identically $0$, (\ref{eq.gradb}) has often appeared \cite{shenm, ks, benali, mp} alongside a strong condition on the size of $V$:
\begin{equation}\label{eq.V}
	V(x)\lesssim m(x,|\B|+V)^2.
\end{equation}
Observe that under assumption (\ref{eq.rh}) for $|\B|+V$, (\ref{eq.gradb}) implies in particular that the components of $\B$ are locally Lipschitz continuous. It is more natural \cite{shenm, benali} to assume (\ref{eq.gradb}) when one is interested in $L^p$ estimates of associated Riesz transforms, for large $p$, since one often needs pointwise bounds for gradients of solutions. However, it would be desirable to dispose of  (\ref{eq.gradb}) when studying questions about the spectrum of $L_{\a,V}$.   Shen remarked \cite[Remark 0.17]{sheno} that for $n\geq3$, some condition on $\nabla\B$ seemed to be necessary to prove a  bound on the number of negative eigenvalues of a magnetic Schr\"odinger operator with a potential $V$ with non-trivial negative part,  but that (\ref{eq.gradb}) was more restrictive than one would hope. 

\begin{conjecture}[Shen's conjecture \cite{sheno}]\label{conj.s1} For $n\geq3$, some condition on $\nabla\B$ is necessary to prove Theorem 0.11 of \cite{sheno}.
\end{conjecture}

However, in the case of $n=2$, Z. Shen \cite{shenb} studied lower and upper bounds for the dimension of the spectral projection of $L_{\a,V}$ on an interval $[0,\tau)$, $\tau>0$,\footnote{The magnetic Schr\"odinger operator on $\bb R^n$ only makes sense for $n\geq2$, due to the gauge invariance.}  under the milder conditions that $\B$ (in this two-dimensional case, given by a single function $b_{12}$) verifies a reverse H\"older inequality and that 
\begin{equation}\label{eq.shensign}
b_{12}  \text{ \emph{ does not change sign} over }\bb R^2.\footnotemark
\end{equation}
\footnotetext{The assumptions on $\B$ in \cite{shenb} are quite more general than stated here.}The question of   extending the ideas of \cite{shenb} to study the analogous problem in higher dimensions   was asked by Shen in that same article \cite[p. 486]{shenb}, but it has remained  open. Let us give a precise statement of the problem that we will consider.

\begin{problem}[Shen's  problem\cite{shenb}]\label{pro.s2} For $n\geq3$, recover the conclusion of \cite[Main Theorem]{shenb} without assuming (\ref{eq.gradb}); instead, impose an assumption on $\B$ that extends, in some sense, the conditions of \cite{shenb} in the two-dimensional setting. 
\end{problem}

The question of relaxing the scale-invariant condition  (\ref{eq.gradb})  when studying the eigenvalue counting function of $L_{\a,V}$ was also asked by S. Mayboroda\footnote{Private communication.},  and is of importance in the case of irregular magnetic fields. In Section \ref{sec.main2a}, we will use a new landscape function for $L_{\a,V}$ to show that, for $n\geq3$, and under favored directionality assumptions (Section \ref{sec.main2}) that in some sense generalize the sign condition (\ref{eq.shensign}), \emph{no assumption on $\nabla\B$ is necessary} (not even merely its existence as a locally integrable function) to study several properties of magnetic Schr\"odinger operators, and in particular, recover the conclusions of \cite[Theorem 0.11]{sheno} and \cite[Main Theorem]{shenb}. We will provide an explicit counterexample (Example \ref{ex.1}) to  Conjecture \ref{conj.s1}, and we will give a solution to Problem \ref{pro.s2}  (see Section \ref{sec.cor}). Moreover, we will not assume \emph{any} quantitative regularity of the magnetic field (not even   $|\B|\in RH_{n/2}$) when stating our a priori Agmon-type results on the exponential decay of eigenfunctions and Lax-Milgram solutions    for $L_{\a,V}$, so that these results hold true also for magnetic fields exhibiting random features.

Let us circle back to the Filoche-Mayboroda landscape function $u$, which solves the equation $-\dv A\nabla u+Vu=1$. As in the case of the maximal function $m$, an uncertainty principle for $u$ akin to (\ref{eq.fpu}) may also be devised, which is in fact quite stronger. In \cite{adfjm}, the following integral identity was shown 
\begin{equation}\label{eq.conj-intro}
	\int_\Omega\Big[A\nabla f\nabla f+Vf^2\Big]=\int_\Omega\Big[u^2A\nabla(f/u)\nabla(f/u)+\frac1uf^2\Big],\qquad f\in C_c^{\infty}(\Omega),
\end{equation}
for a bounded domain $\Omega\subset\bb R^n$, $n\geq3$, with Neumann boundary conditions, symmetric $A$, and $V\in L^{\infty}(\Omega)$. Observe that (\ref{eq.conj-intro}) readily implies the uncertainty principle
\begin{equation}\label{eq.up0}
\int_{\Omega}\frac1uf^2\leq\int_{\Omega}\Big[A\nabla f\nabla f+Vf^2\Big],
\end{equation}
which ought to be compared to (\ref{eq.fpu}). As such, it is not surprising that (\ref{eq.up0}) was used along with an Agmon argument to give the exponential decay of eigenfunctions to $L$ on $\Omega$, with respect to the Agmon distance with weight $1/u$. On the other hand, if $V$ also verified the reverse H\"older inequality (\ref{eq.rh}), then we may also  give the exponential decay of eigenfunctions with an Agmon distance with weight $m^2(\cdot,V)$. It is thus natural to wonder, if there is any relationship between the landscape function $u$ and the   maximal function $m$, at least in the cases when the latter makes sense. A main result of our paper, which is at the heart of this investigation, is that  such a relationship does exist (see Theorem \ref{thm.compare}).

\subsection{Main results for the Schr\"odinger operator  $-\dv A\nabla+V$}\label{sec.main} Given an open, connected set $\Omega\subset\bb R^n$, $n\geq1$, with (possibly empty) Lipschitz boundary, we let $0\leq V\in L^1_{\loc}(\Omega)$ with $\int_\Omega V>0$, and $G=G_{\Omega}$ be the Green's function for the operator $L$ on  $\Omega$ with $0$ Dirichlet boundary conditions (see Section \ref{sec.prelim}). Under our very mild assumptions on $V$, the Green's function may not exist as a distributional solution to the equation $LG(x,\cdot)=\delta_x$ in $\Omega$ \cite{bb03, pw17}, but it always exists as a non-negative\footnote{The non-negativity of $G$ follows essentially from the non-degeneracy condition $\int_\Omega V>0$.} integral kernel of the (homogeneous) operator $L^{-1}$ (see   Appendix \ref{sec.green}).

When $\Omega$ is bounded, the landscape function may be easily taken as the Lax-Milgram solution $L^{-1}{\1}_{\Omega}$ (see Section \ref{sec.prelim} for the definition of $L^{-1}$), but this approach does not work on unbounded domains. We are interested in the unbounded case because later it will allow us to connect the landscape function to the maximal function $m$ on $\Omega=\bb R^n$, which is critical to our main results for $L_{\a,V}$.  The first main result  of our paper is to establish that the existence of a landscape function over unbounded domains is equivalent to a mild decay assumption on the Green's function. If this decay condition holds, then we show that, even though in our setting the zero-set of $u$ in the interior of $\Omega$ may not be trivial, and even though it may not be true that $\frac1u\in L^{\infty}_{\loc}(\Omega)$, nevertheless, the uncertainty principle (\ref{eq.up0}) remains valid.

\begin{theorem}[A landscape function over unbounded domains]\label{thm.exist}  Let $n\geq1$, let $\Omega\subseteq\bb R^n$ be an open   connected set, with empty or Lipschitz $(n-1)-$dimensional boundary,  let $A$ be a not necessarily symmetric matrix of bounded measurable coefficients on $\Omega$ verifying the uniform ellipticity condition (\ref{eq.elliptic}), and   $0\leq V\in L^1_{\loc}(\overline{\Omega})$ satisfies $\int_\Omega V>0$. Denote $L=-\dv A\nabla+V$. The following statements hold.
\begin{enumerate}[(i)]
\item\label{item.u} Fix $x_0\in\Omega$ and for each $R\in\bb N$, let $\Omega_R:=\Omega\cap B(x_0,R)$. Consider the landscape function $u_R:=L^{-1}_{\Omega_R}\1_{\Omega_R}$ on the domain $\Omega_R$  with zero Dirichlet boundary conditions, and extend it by $0$ on $\Omega\backslash\Omega_R$. Then the sequence $\{u_R\}_{R=1}^{\infty}$ is pointwise non-decreasing, and the limit
\begin{equation}\label{eq.limit}
u:=\limsup_{R\ra\infty}u_R
\end{equation}
exists as a measurable non-negative function on $\Omega$ (whose values may be $+\infty$ everywhere).
\item\label{item.equiv} The function $u$ of \ref{item.u} verifies that $u\in W^{1,2}_{\loc}(\Omega)\cap L^2_{\loc}(\Omega,V\,dx)$ and solves the equation $Lu=1$ in the weak sense in $\Omega$ (see Definition \ref{def.weak}) if and only if there exists $q>0$ such that
\begin{equation}\label{eq.greenint}
\int_{\Omega}G_{\Omega}(x,y)\,dy\in L^q_{\loc}(\Omega,\,dx).
\end{equation}
\item\label{item.under} If (\ref{eq.greenint}) holds, then we call $u$ the \emph{landscape function} for $L$ on $\Omega$ (with zero Dirichlet boundary conditions), and moreover, we have that $u=\int_\Omega G(\cdot,y)\,dy$ a.e. on $\Omega$, and the uncertainty principle (\ref{eq.unprinciple}) holds for each $f\in \f D(\Omega)$\footnote{See Section \ref{sec.prelim} for the definition of the space $\f D(\Omega)$.}. In particular, $\frac1u\in L^1_{\loc}(\Omega)$ and $\nabla\log(u)\in L^2_{\loc}(\Omega)$. 

Furthermore, if $A$ is symmetric, then the stronger uncertainty principle (\ref{eq.unprinciplesym}) holds for each $f\in \f D(\Omega)$.
\end{enumerate}
\end{theorem}

Of course, if $\Omega$ is bounded, then (\ref{eq.greenint}) holds by  (\ref{eq.greenlp}). In the case of unbounded $\Omega$, the estimate (\ref{eq.greenint}) is a very mild decay condition on the Green's function, which implies some non-degeneracy of $V$ at infinity. Indeed, if $V=0$ outside of a bounded set in $\bb R^n$ and $V\in L^{\infty}(\bb R^n)$, then (\ref{eq.greenint}) is not satisfied, and the limit in (\ref{eq.limit}) is $+\infty$ everywhere on $\bb R^n$. On the other hand, for $n\geq3$, $\Omega=\bb R^n$, and for $V$ verifying (V1) and (V2) of Theorem \ref{thm.compare} below, then (\ref{eq.greenint}) is known to hold, by virtue of the exponential decay results of the fundamental solution \cite{shenf, mp}\footnote{The results in \cite{mp} are shown under stronger hypotheses, but it is an easy exercise to see that they extend to the Shen potentials, following \cite{shenf}.}.

In fact, when $A\equiv I$, the assumptions on $V$ allow for the existence of a \emph{universal zero-set} $Z\subset\Omega$ such that every weak solution of $Lu=f$, $f\in L^{\infty}(\Omega)$ vanishes on $Z$ \cite[Examples 1.1-1.4]{op2020}, so that the strong maximum principle does not hold in our generality. In particular, the landscape function that we construct may not be continuous, and is forced to take the value $0$ at points in the interior of $\Omega$ where $V$ is highly singular \cite[Theorem 1.1 and Corollary 1.2]{op2020}; both of these properties differ substantially from the setting of \cite{adfjm}. Still,   motivations to consider such singular $V$ include that, by the weak maximum principle we should expect \emph{more} exponential decay of solutions when $V$ is higher, and we should expect that the landscape function can still quantify this decay; and also, $L^1_{\loc}$ is the minimal condition that can be imposed on a measurable function so that the bilinear form induced by $L$ acts on smooth compactly supported functions.

We remark that the continuity of the landscape function is always true in dimension $1$ by Morrey's inequality, and the strong maximum principle can be verified for $n\geq1$ if $V$ belongs to the Lorentz space $L^{\frac n2,1}_{\loc}(\Omega)$ or if $V$ belongs to the Kato class $\m K_n$. The issue of whether the strong maximum principle holds is strongly tied to certain capacities of the zero sets  of solutions \cite{ancona79}. For more on this, see \cite{op16, op2020}.

With the existence of a landscape function at hand, we tackle \emph{a priori}  exponential decay estimates for Lax-Milgram solutions,   Green's functions, and eigenfunctions, extending  results of \cite{adfjm} to the case $V\in L^1_{\loc}(\Omega)$ and to Dirichlet boundary conditions\footnote{See also Corollary \ref{cor.resexp} for the exponential decay of resolvents.}. The a priori exponential decay of eigenfunctions for tight-binding Hamiltonians has also been considered \cite{wz21}. The definitions of   Agmon distances $\rho,\rho_A$ and the  version  $\hat u$ of the landscape function are found in Section \ref{sec.exp}. 

\begin{theorem}[Exponential decay estimates I: Lax-Milgram solutions]\label{thm.lmr} Retain the setting of Theorem \ref{thm.exist},   assume that (\ref{eq.greenint}) holds, and that
\begin{equation}\label{eq.nosing}
\rho\big(x,y,\frac1{\hat u}\big)\lra0,\quad\text{as }|x-y|\ra0,\qquad\text{for any }x,y\in\Omega.
\end{equation}
Then there exist constants $\ep, C$, depending only on $\lambda$, such that for every $f\in L^2(\Omega)$ with compact support in $\Omega$, 
\begin{equation}\label{eq.lmest}
\int_{\Omega}\frac1{u}e^{2\ep\rho(\cdot,\supp f,\frac1{\hat u})}|L^{-1}f|^2~\leq C\int_{\Omega}uf^2.
\end{equation}
Moreover, if $A$ is symmetric and continuous, then for each $\alpha\in(0,1/4)$,
\begin{multline}\label{eq.lmests}
	\int_{\Omega}u^2A\nabla\Big(\frac{e^{\alpha\rho_A(\cdot,\supp f,\frac1{\hat u})}L^{-1}f}{u}\Big)\nabla\Big(\frac{e^{\alpha\rho_A(\cdot,\supp f,\frac1{\hat u})}L^{-1}f}{u}\Big)\\+\int_{\Omega}\frac1{u}e^{2\alpha\rho_A(\cdot,\supp f,\frac1{\hat u})}|L^{-1}f|^2~\leq\frac1{(1-\alpha^2)^2}\int_{\Omega}uf^2.
\end{multline} 
\end{theorem}

Assumption (\ref{eq.nosing}) says that the rate at which $u\searrow0$ within compact subsets of $\Omega$ is mildly controlled; for instance, this assumption holds if $u>0$ on $\Omega$ (and this is true, in particular, for bounded potentials, or potentials in the Kato class $\m K_n$).

Knowing that the Green's function of an operator exhibits large decay is   useful in the study of the spectrum of the operator, as well as towards proving $L^p$ estimates of associated Riesz transforms, and related questions in harmonic analysis \cite{shenm, benali, bailey2021, bmr21}. The following theorem says that the decay of the integral kernel is a priori governed by the behavior of the landscape function.

\begin{theorem}[Exponential decay estimates II: Green's function]\label{thm.greenexp} Retain the setting of Theorem \ref{thm.exist}, and assume that (\ref{eq.greenint}) and (\ref{eq.nosing}) hold. Then, for each $\alpha\in(0,1)$, we have that
\begin{multline}\label{eq.gineqavg}
\dashint_{B(x_0,\delta)}\dashint_{B(y_0,\ep)}G(x,y)\,dy\,dx\\ \leq\frac1{|B(0,1)|}\frac1{(1-\alpha^2)}\delta^{-\frac n2}\ep^{-\frac n2}\Big(\dashint_{B(y_0,\ep)}u\Big)^{\frac12}\Big(\dashint_{B(x_0,\delta)}u\Big)^{\frac12}\max_{x\in B(x_0,\delta)}\min_{y\in B(y_0,\ep)}e^{-\alpha\rho(x,y,\frac1{\hat u})}
\end{multline}
for each $x_0,y_0\in\Omega$ and $\delta,\ep>0$ such that $B(x_0,\delta)\cup B(y_0,\ep)\subset\Omega$. 

Furthermore, if $n\geq3$ and $V\in L^{\infty}(\Omega)$, then for any $x,y\in\Omega$ with
\[
B_x\cap B_y:=B\big(x, \Vert V\Vert_{L^{\infty}(\Omega)}^{-1/2}\big)\cap B\big(y,\Vert V\Vert_{L^{\infty}(\Omega)}^{-1/2}\big)=\varnothing,\qquad 4B_x\cup4B_y\subset\Omega,
\]
there exist $\tilde x\in B_x$, $\tilde y\in B_y$ such that the estimate
\begin{equation}\label{eq.greendecay}
G(x,y)\leq C\Vert V\Vert_{L^{\infty}(\Omega)}^{\frac n2}\sqrt{u(\tilde x)}\sqrt{u(\tilde y)}e^{-\alpha\rho(\tilde x,\tilde y,\frac1{u})}
\end{equation}
holds, and  $C$ depends only on $n$, $\alpha$, and $\lambda$.
\end{theorem}

For the a priori exponential decay estimate on eigenfunctions, we restrict to symmetric and continuous $A$. The existence of a point spectrum $\sigma_p$ for  Schr\"odinger operators, and the conditions on $V$ or $\Omega$ which guarantee it, are  well-researched topics in the literature. However, due to the generality of our setting, in this paper we will not deal with the  existence of $\sigma_p$, and we will simply show that if an eigenfunction exists, then its decay is governed by the Agmon distance with  a weight given by the landscape function.

\begin{theorem}[Exponential decay estimates III: Eigenfunctions]\label{thm.eigen} Retain the setting of Theorem \ref{thm.exist}, and moreover, assume that (\ref{eq.greenint}) and (\ref{eq.nosing}) hold, and that $A$ is symmetric and continuous\footnote{If $A$ is not continuous, then the results here still hold with $\rho_A$ replaced by $\rho$, and up to multiplicative constants depending on $\lambda$.}. Let $\n L$ be the operator   in Definition \ref{def.nonhomo}, with domain $\n D(\n L)\subset L^2(\Omega)$. Suppose that there exist $\mu>0$ and $\psi\in\n D(\n L)$ such that $\n L\psi=\mu\psi$. Let
\[
w(x):=\big(\tfrac1{\hat u(x)}-\mu\big)_+=\max\big\{0,\tfrac1{\hat u(x)}-\mu\big\},\quad\text{and}\quad E:=\big\{x\in\Omega: \tfrac1{\hat u(x)}\leq\mu\big\}.
\]
Then for each $\alpha\in(0,1/4)$, we have that
\begin{multline}\label{eq.eigen}
\int_{\Omega}u^2A\nabla\Big(\frac{e^{\alpha\rho_A(\cdot,E,w)}\psi}{u}\Big)\nabla\Big(\frac{e^{\alpha\rho_A(\cdot,E,w)}\psi}{u}\Big)\\+\int_{\Omega}\Big(\frac1{u}-\mu\Big)_+e^{2\alpha\rho_A(\cdot,E,w)}|\psi|^2~\leq\frac1{(1-\alpha^2)}\int_E\Big(\mu-\frac1u\Big)_+|\psi|^2.
\end{multline} 
\end{theorem}

\begin{remark}  Note that, under the assumptions of Theorem \ref{thm.eigen}, the right-hand side of (\ref{eq.eigen}) is always finite, since $\psi\in L^2(\Omega)$ and $\int_\Omega\frac1u\psi^2<+\infty$ because of the uncertainty principle (\ref{eq.unprinciplesym}). 
\end{remark}

The proofs of the three previous theorems may be found in Section \ref{sec.exp}  and consist  of  an Agmon-type argument, using the uncertainty principle, similar as in \cite{adfjm}. However,  technical issues have to be dealt with since $V$ is no longer   bounded (except for (\ref{eq.greendecay})), $\frac1u\ra+\infty$ at $\partial\Omega$, and   $u$ is not even necessarily positive on $\Omega$. The averaged $L^1$ decay result of Theorem \ref{thm.greenexp}   relies on Theorem \ref{thm.lmr} and the representation formula for Lax-Milgram solutions in terms of the Green's function. From (\ref{eq.gineqavg}) and a lower bound for $u$ in terms of $V$ (which is not trivial when $u$ is the Dirichlet landscape function; see Proposition \ref{prop.harnacku}), (\ref{eq.greendecay}) follows. Let us emphasize that the decay estimates of the previous theorems do not assume any regularity of the potential $V$; hence, in particular, they hold for random potentials. On the other hand, in the generality of our setting  the Agmon distance $\rho(x,y,\frac1u)$ may be degenerate, and may not look anything like the Euclidean distance, even for bounded  $V$. This is why we say that our estimates are a priori.

Now we restrict to $n\geq3$ for the next few results, which will give situations where the Agmon distance with weight $1/u$ is not degenerate. Of course, the estimates of Theorems \ref{thm.lmr} and \ref{thm.greenexp} also hold for potentials $V\in RH_{n/2}$, and in this case these results overlap with certain upper bound exponential decay estimates from \cite{shenf, mp}, where the maximal function $m$ was used; however, the Agmon distance used in these papers a priori differs from the Agmon distance used in this article, because $m^2$ and $1/u$ are different weights. Nevertheless, this is not an issue: our next main result  shows that the functions $m^2$ and $1/u$ are pointwise equivalent. In fact, we show this under the much milder assumptions of \cite{shenf}, themselves inspired by \cite{christ91}.

\begin{theorem}[The relation between $1/u$ and $m^2(\cdot,V)$]\label{thm.compare} Retain the setting of Theorem \ref{thm.exist}, and moreover, assume that $n\geq3$, that $\Omega=\bb R^n$, and  suppose that $V$ verifies:
	\begin{enumerate}[(V1)]
		\item (Scale-invariant Kato Condition). There exist positive constants $C_0$ and $\delta$ such that for all $x\in\bb R^n$ and all $0<r<R$,
		\begin{equation}\label{eq.kato}
			\int_{B(x,r)} V(y)\,dy\leq C_0\Big(\frac rR\Big)^{n-2+\delta}\int_{B(x,R)}V(y)\,dy.
		\end{equation}
		\item (Doubling on balls with high mass). There exists a positive constant $C_1$ such that for all $x\in\bb R^n$ and all $r>0$,
		\begin{equation}\label{eq.largedoubling}
			\int_{B(x,2r)}V(y)\,dy\leq C_1\Big[\int_{B(x,r)}V(y)\,dy+r^{n-2}\Big].
		\end{equation}
	\end{enumerate}

Then (\ref{eq.greenint}) holds, and there exists   $C>0$, depending only on $n,\lambda, C_0, C_1, \delta$, such that 
	\begin{equation}\label{eq.compare}
		\frac1C\frac1{m^2(x,V)}\leq u(x)\leq C\frac1{m^2(x,V)},\qquad\text{for each }x\in\bb R^n,
	\end{equation}	
where $m(\cdot,V)$ is the maximal function from (\ref{eq.fpsm}), and $u$ is the landscape function defined in Theorem \ref{thm.exist}.
\end{theorem}

\begin{remark} It is important in Theorem \ref{thm.compare} that we take the solution $u$ to $Lu=1$ described in Theorem \ref{thm.exist}; other solutions to $Lv=1$ may not verify (\ref{eq.compare}), even though they still verify the uncertainty principle (\ref{eq.unprinciple}) for $C_c^{\infty}(\bb R^n)$ test functions. For instance, take $V\equiv1$; in this case, $m(\cdot,1)=C$, and consider $v=e^{x_1}+1$, which solves $-\Delta v+v=1$ on $\bb R^n$. It is clear that (\ref{eq.compare}) does not hold for $v$ in place of $u$.	
\end{remark}

The potentials verifying (V1) and (V2) will be called \emph{Shen potentials}. The proof of this theorem is given in Section \ref{sec.fps}, and relies on exponential decay estimates from \cite{mp} for Lax-Milgram solutions. 

Let us make a few remarks regarding the conditions (\ref{eq.kato}) and (\ref{eq.largedoubling}). First, if $V$ verifies the reverse H\"older inequality (\ref{eq.rh}), then $V$ is a Shen potential. On the other hand,  $V(x)=e^{-|x|}$ is a Shen potential, but does not lie in any reverse H\"older class $RH_q$, $q>1$ (because $V$ is not doubling). More generally, note that if $V\in RH_q$, then $V$ is the Radon-Nikodym derivative of a measure which is $A_{\infty}-$absolutely continuous with respect to the Lebesgue measure, but Shen potentials may be zero on open sets, for instance. Finally, let us note that random Schr\"odinger potentials fail both conditions (V1) and (V2).

There is a practical advantage in working with $1/u$ instead of $m^2$, even though they are pointwise equivalent by the previous theorem. Namely, from the point of view of numerics, it is much easier to compute the solution to the problem $Lu=1$ than to compute the maximal function $m$; moreover, under the setting of Theorem \ref{thm.compare}, $u$ is H\"older continuous, while $m$ is not. Now, from   well-known results concerning  $m$, we are able to give   immediate corollaries which are new for the landscape function. The first one uses the exponential decay results of \cite{shenf, mp} to obtain upper and lower bound exponential decay estimates in terms of an Agmon distance with weight $1/u$.

\begin{corollary}[Sharp exponential decay with the landscape function for Shen potentials]\label{cor.sharp} Retain the setting and assumptions of Theorem \ref{thm.compare}, and let $\Gamma$ be the fundamental solution for the operator $L$. Then there exist constants $C$, $\ep_1$, and $\ep_2$, which depend only on $n$, $\lambda$, $C_0$, $C_1$, and $\delta$, such that
	\begin{equation}\label{eq.sharp}
		\frac1{C|x-y|^{n-2}}e^{-\ep_1\rho(x,y,\frac1u)}\leq\Gamma(x,y)\leq C\frac1{|x-y|^{n-2}}e^{-\ep_2\rho(x,y,\frac1u)},\qquad\text{for each }x,y\in\bb R^n.
	\end{equation}
\end{corollary}

The next corollary gives  properties that $u$ enjoys by virtue of being pointwise equivalent to $1/m^2$, under the assumptions (\ref{eq.kato}) and (\ref{eq.largedoubling}).

\begin{corollary}[Properties of the landscape function for Shen potentials]\label{cor.u} Retain the setting and assumptions of Theorem \ref{thm.compare}. Then there exist   $C>0$ and $k_0\in\bb N$, which depend only on $n$, $\lambda$, $C_0$, $C_1$, and $\delta$, such that the following statements hold.
\begin{enumerate}[(i)]
\item \emph{($1/u$ encodes $V$).} For each $x\in\bb R^n$, we have that
\begin{equation}\label{eq.comparetoV}
	\frac1C\frac1{u(x)}\leq\dashint_{B(x,\sqrt{u(x)})}V(y)\,dy\leq C\frac1{u(x)}
\end{equation}
\item \emph{(Scale-invariant Harnack inequality)}. For each $x\in\bb R^n$, if $y\in B(x,\sqrt{u(x)})$, then
\begin{equation}\label{eq.harnacku}
	\frac1Cu(x)\leq u(y)\leq Cu(x).
\end{equation}
\item \emph{(Comparability at a long distance).} For each $x,y\in\bb R^n$, we have that
\begin{equation}\label{eq.long}
	\frac1C\frac1{\sqrt{u(y)}} \Big\{1+\frac{|x-y|}{\sqrt{u(y)}}\Big\}^{\frac{-k_0}{k_0+1}} \leq\frac1{\sqrt{u(x)}}\leq C\frac1{\sqrt{u(y)}}\Big\{1+\frac{|x-y|}{\sqrt{u(y)}}\Big\}^{k_0}.
\end{equation}
\item \emph{(Case of polynomial potentials \cite[Estimate (0.8)]{sheno}).} If $V(x)=|P(x)|^{\alpha}$ where $P$ is a polynomial of degree $k\in\bb N$ and $\alpha>0$, then
\begin{equation}\label{eq.polyu}
\frac1{u(x)}\approx\Big(\sum_{|\beta|\leq k}|\partial_x^{\beta}P(x)|^{\frac{\alpha}{\alpha|\beta|+2}}\Big)^2,\qquad\text{for each }x\in\bb R^n.
\end{equation}
\end{enumerate} 
\end{corollary}

The scale-invariant Harnack inequality (\ref{eq.harnacku}) is also known as the \emph{slowly varying} property; note that this property is significantly stronger than what the regular Harnack inequality \cite[Theorem 8.18]{gtpde} would predict for $u$. The estimates (\ref{eq.long}) imply that the Agmon distance with respect to $1/u$ is pointwise controlled by powers of the Euclidean distance. Therefore, under the assumptions of Theorem \ref{thm.compare}, exponential decay in the Agmon distance with weight $1/u$ \emph{also implies exponential decay in the Euclidean distance}.   To the author's knowledge, this is the first time that the link between exponential decay in the Agmon distance with weight $1/u$ and the exponential decay in the Euclidean sense has been rigorously investigated. On the other hand, let us emphasize that we have only proved the results in the last corollary for potentials which verify the assumptions (\ref{eq.kato}) and (\ref{eq.largedoubling}), and random potentials decidedly do not satisfy either of these assumptions.

\subsection{Directionality assumptions on the magnetic field}\label{sec.main2} Let us now show an application of the previous theorems.  We construct  a new landscape function  for the magnetic Schr\"odinger operator $L_{\a,V}$ given in  (\ref{eq.ms}) on an open domain $\Omega\subset\bb R^n$, $n\geq2$, and prove that it gives exponential decay of solutions\footnote{More precisely, we construct a candidate for a landscape function; the question of whether the function constructed here is helpful in the study of Anderson localization  under the presence of a random magnetic field will be the subject of future investigations.}. Here, the real-valued vector function ${\bf a}(x)=(a_1(x),\ldots,a_n(x))$ is the \emph{magnetic potential}, and recall that the magnetic field $\B$ is given in (\ref{eq.magfield}). 

Our idea in constructing a landscape function for $L_{\a,V}$ relies on the \emph{gauge invariance} property of this operator: if $\Phi\in C^1$ is any real-valued function, then the magnetic potential $\tilde{\bf a}={\bf a}+\nabla\Phi$ has the same magnetic field as that of ${\bf a}$; hence $\tilde{\bf B}={\bf B}$. Moreover, the identities
\begin{equation}\label{eq.gauge}
	L_{\tilde{\bf a}}f=e^{-i\Phi}L_{\bf a}(e^{i\Phi}f),\qquad D_{j,\tilde{\bf a}}f=e^{-i\Phi}D_{j,{\bf a}}(e^{i\Phi}f)
\end{equation}
hold, where
\begin{equation}\label{eq.maggrad}
D_{j,\bf a}f:=\frac{\partial f}{\partial x_j}-ia_jf,\qquad 1\leq j\leq n.
\end{equation}
The energy associated to the operator (\ref{eq.ms}) is
\begin{equation}\label{eq.energy}
	\int_{\bb R^n}\Big[|D_{\bf a}f|^2+V|f|^2\Big],\qquad f\in C_c^{\infty}(\bb R^n),
\end{equation}
where
\[
D_{\bf a}:=(D_{j,{\bf a}})_{1\leq j\leq n},
\]
and thus observe that, by (\ref{eq.gauge}), the energy associated to $L_{\bf a}$ is invariant under gauge transformations, and determined solely by the magnetic field ${\bf B}$. For this reason, when seeking a landscape function for the magnetic Schr\"odinger operator that generalizes the Filoche-Mayboroda landscape function for the non-magnetic case, it makes sense to attempt to create it solely by knowledge of the magnetic field, and not of the specific magnetic potential ${\bf a}$, which can always be ``gauged out''. Moreover, we seek to do this under minimal assumptions on the magnetic field; in particular, we do not want to assume any pointwise bound like (\ref{eq.gradb}) on the gradient of the magnetic field, since we want to be able to treat  random magnetic potentials as well (at least, in an a priori  way, similar as how it was done in Section \ref{sec.main}). Let us also note that (\ref{eq.energy}) is invariant over orthonormal changes of variables.

We now state our assumptions on the magnetic field. Without loss of generality we may assume that the origin lies in $\Omega$. We say that $\m S=\{(j',k')\}$ is a \emph{selection of pairs} if for each $j,k\in\{1,2,\ldots,n\}$ with $j \neq k$, exactly one of $(j,k)$ or $(k,j)$ belong to $\m S$. Selections of pairs will correspond to choosing exactly one function of $b_{jk}$ or $b_{kj}$, for each $j,k$. Note that, since $\B$ is antisymmetric,
\[
|\B|(x)=\sum_{(j',k')\in\m S}|b_{j'k'}(x)|,\quad\text{for any selection of pairs }\m S.
\]
We assume that
\begin{equation}\label{eq.bl1}
  {\bf a}\in L^2_{\loc}(\Omega), b_{jk},V\in L^1_{\loc}(\Omega),\quad\text{for each }j,k=\{1,2,\dots,n\},
\end{equation}
and that
\begin{multline}\label{eq.signass}
\text{possibly after an orthonormal change of coordinates, there exists } \m S \text{ such that }\\ \sum_{(j',k')\in\m S}b_{j'k'}+V\geq0 \text{ pointwise a.e. on }\Omega.
\end{multline}
If $\m S$ is such that the inequality in (\ref{eq.signass}) is verified (for some given coordinates), then we call $\m S$ an \emph{admissible selection}.

Observe that, if $V=0$, (\ref{eq.signass}) is essentially  a requirement that the magnetic field  favors certain directions.  For instance,  (\ref{eq.signass}) is verified if   
\begin{multline}\label{eq.ass}
|\B|+V\geq0  \text{ and, possibly after an orthonormal change of coordinates,}\\ b_{jk}  \text{ does not change sign over }\Omega, \text{ for each }j,k=\{1,2,\ldots,n\}.
\end{multline}
In the case that  (\ref{eq.ass}) holds, there is a unique admissible selection $\m S$ which verifies
\begin{equation}\label{eq.magb2}
	|\B|=\sum_{(j',k')\in\m S} b_{j'k'}.
\end{equation} 

We refer to (\ref{eq.signass}) and (\ref{eq.ass}) respectively as the weak and strong \emph{favored directionality assumptions}. Let us discuss the meaning and viability of these conditions.

For $n=2$, (\ref{eq.ass}) is easily seen to be related to  (\ref{eq.shensign}) which appeared in \cite{sheno, bls2004}, and if $V=0$, it simply states that the direction of the magnetic field is not allowed to flip; (\ref{eq.shensign}) was also considered in \cite{eh2012} for the study of Anderson localization for random magnetic fields, and the Aharonov-Bohm fields \cite{mn2009} also formally\footnote{For the Aharonov-Bohm fields, the magnetic field is a measure with atoms.} verify the assumption. For $n\geq3$, the landmark paper \cite{ahs78} introduced a condition on the direction of the magnetic field in their Corollary 2.10 which allowed the construction of magnetic bottles over $\bb R^n$. They gave as an example  that in three dimensions where we may identify the magnetic field as $(b_1,b_2,b_3)$, if $b_1,b_2,b_3$ are non-negative and $b_1+b_2+b_3\ra+\infty$ as $|x|\ra\infty$, then the operator $L_{\a,0}$ has compact resolvent. Thus the favored directionality assumptions (\ref{eq.signass}), (\ref{eq.ass}) are very much related to the condition of \cite[Corollary 2.10]{ahs78}. This latter condition  has been used many times in the literature \cite{v86, es99, es04, es04b, truc97, vt2010, nakamura2000, bls2004} to obtain uncertainty principles, magnetic Hardy inequalities, and magnetic bottles; on the other hand, these papers have often introduced extra \emph{regularity} assumptions  on the direction of the magnetic field, requiring its smoothness in a scale-invariant sense, which is what we are trying to avoid in this paper (but see Theorem \ref{thm.transition} for a case where we consider a mixture of our assumption (\ref{eq.ass}) and the regularity condition (\ref{eq.gradb})).

Let us give two intuitive reformulations of (\ref{eq.ass}) for the physically relevant case $n=3$. Here, the magnetic field may be identified by a vector function ${\bf b}=(b_1,b_2,b_3)$, and in this case the second statement of (\ref{eq.ass}) is equivalent to  
\begin{align}\nonumber
\text{the range of }{\bf b}\text{ is a subset of the convex hull}& \text{ of the positive half-axes}  \\ &\text{ associated to some orthonormal basis of }\bb R^3.\nonumber
\end{align}
Another equivalent way to understand the second statement of (\ref{eq.ass}) for $n=3$ is
\[
{\bf b}(x)\cdot{\bf b}(y)\geq0,\qquad\text{for each }x,y\in\Omega.
\]
Therefore, (\ref{eq.ass}) says that  the direction of the magnetic field must be somewhat consistent over the whole domain. Note that for any $n$, if $\B$ is continuous and if $V\geq0$, then (\ref{eq.ass}) is always satisfied in a neighborhood of a point where the multiplicity of $0$ as an eigenvalue for $\B$ is no more than $1$; this is a consequence of the spectral theory for antisymmetric matrices. Let us be more precise: fix $x_0\in\Omega$, and since $\B(x_0)$ is a real-valued antisymmetric matrix, there is an orthonormal change of coordinates $\Theta$ (depending on $x_0$) so that, if $n=2m$ is even, 
\[
\B^{\Theta}(x_0)=\begin{bmatrix}
	\begin{matrix}0 & b_{12}^{\Theta} \\ -b_{12}^{\Theta} & 0\end{matrix} &  0 & \cdots & 0 \\
	0 & \begin{matrix}0 & b_{34}^{\Theta}\\ -b_{34}^{\Theta} & 0\end{matrix} & & 0 \\
	\vdots & & \ddots & \vdots \\
	0 & 0 & \cdots & \begin{matrix}0 & b_{(m-1)m}^{\Theta}\\ b_{(m-1)m}^{\Theta} & 0\end{matrix}
\end{bmatrix}
\]
with $b_{12}^\Theta\geq b_{34}^\Theta\geq\cdots\geq b_{(m-1)m}^\Theta\geq0$. If $n$ is odd, we have to add a row and column of $0$ to $\B^\Theta(x_0)$, but we still get the same sequence $\{b_{j(j+1)}^\Theta\}$ as above, with $n=2m+1$. It is clear then that (\ref{eq.ass}) is verified at $x_0$, and if $b_{(m-1)m}>0$, then by the continuity of $\B$ we get that (\ref{eq.ass}) is verified in a small neighborhood of $x_0$.  

We mention that, for $n\geq2$, the assumptions (\ref{eq.signass}) and (\ref{eq.ass}) are well-suited to work with random Landau Hamiltonians, where the magnetic field is constant. They are also suited to work with random magnetic potentials of the form $\a_0+\delta\a_{\omega}$, where $\a_0$ is a given magnetic potential verifying (\ref{eq.signass}), $\a_{\omega}$ is an arbitrary bounded random magnetic potential, and $\delta>0$ is sufficiently small. These situations have been considered heavily in the literature; some related works are \cite{ch96, hlmw01, kr06, ghk2007, klopp10}.

Given a magnetic field $\B$ over $\bb R^n$, it is  easy to consider many examples when $\Omega$ is chosen so that (\ref{eq.ass}) holds for $\B$ on $\Omega$. Examples of (\ref{eq.ass}) over $\bb R^n$ include magnetic fields whose components are polynomials which do not change signs (in particular, constant magnetic fields verify (\ref{eq.ass})), but these are treated  by existing literature, since they verify the smoothness assumption (\ref{eq.gradb}). Let us now give an  example over $\bb R^n$ which is not handled by the assumptions of Z. Shen \cite{shene,shenb,sheno,mp,kurata}, but which satisfies our assumptions.

\begin{example}\label{ex.1}  For $n\geq3$, let  $\alpha\in(0,1)$, and consider   $\a=(a_1,\ldots,a_n)$ given by
\[
a_j(x)=\left\{\begin{matrix}x_{j+1}^\alpha,\quad &x_{j+1}\geq0&\\-(-x_{j+1})^{\alpha},\quad &x_{j+1}<0,&\end{matrix}\right.\qquad j=1,2,\ldots,n,\quad x_{n+1}\equiv x_1.
\]
Clearly, $\a$ is H\"older continuous, but not Lipschitz continuous, since
\[
\frac{\partial a_j}{\partial x_{j+1}}=\alpha|x_{j+1}|^{\alpha-1},\quad\frac{\partial a_j}{\partial x_k}=0,\qquad j=1,\ldots,n, k\neq j+1, \quad x_{n+1}\equiv x_1.
\]
It follows that $b_{j(j+1)}=-b_{j(j-1)}=\alpha|x_{j+1}|^{\alpha-1}$, and $b_{jk}=0$ otherwise for $j=1,\ldots, n$, and where we have identified the subscripts $n+1$ with $1$, and $0$ with $n$. Let $\m S=\{(j,k):j<k\}$ be the upper triangular selection. Then this example verifies assumption (\ref{eq.ass}), and
\[
|\B|=\alpha\sum_{j=1}^n|x_j|^{\alpha-1}=\sum_{(j',k')\in\m S}b_{j'k'},
\]
so that $\m S$ is the unique admissible selection from (\ref{eq.magb2}). Notice that $|\B|$ is smooth away from the axes, and $|\B|\in L^p_{\loc}(\bb R^n)$ for each $p\in(1,\frac1{1-\alpha})$. A simple calculation\footnote{It is easier to show the property $(\int_Q|\B|^p)^{\frac1p}\lesssim\int_Q|\B|$ over cubes $Q$ in this example.} will show also that $|\B|\in RH_p$ for any such $p$.  In particular, if $\alpha>\frac{n-2}n$, then $|\B|\in RH_{n/2}$. On the other hand,  a weak derivative does not exist as an integrable function near the coordinate axes for any non-zero component of $\B$. 

Finally, note that if $\alpha>\frac{n-2}n$ and $-\frac{\alpha}{2n}\sum_j|x_j|^{\alpha-1}\leq V\leq0$, then we still have $0\leq|\B|+V\in RH_{n/2}$, so that the maximal function $m(\cdot,|\B|+V)$ is well-defined, and this example verifies the assumptions of Theorem  \ref{thm.magup}. If $V\geq0$, then this example furthermore verifies the assumptions of Theorem \ref{thm.expdecayms} (except for the assumptions in \ref{item.greenms}, since the magnetic field is unbounded) and Corollaries \ref{cor.bounds} and \ref{cor.neg}. \hfill{$\square$}
\end{example}

Note that the condition (\ref{eq.signass}) is quite a bit more general than (\ref{eq.ass}), since it allows for  cancellations to occur between the different components of the magnetic field and the potential $V$ as well. In particular, we do not in general require $V\geq0$ in (\ref{eq.signass}) nor in (\ref{eq.ass}), so that these assumptions allow one to study cases when the negative contributions of $V$ are cancelled by the magnetic field, as illustrated in Example \ref{ex.1}. In the opposite case where $V$ is dominant, note that if
\begin{equation}\label{eq.Vdominant}\nonumber
	|\B|\leq \frac2{(n-1)n}V,
\end{equation}
then (\ref{eq.signass}) is verified in any coordinate system, and any selection $\m S$ is an admissible selection\footnote{This situation, where the magnetic field is dominated by $V$, was studied in \cite{benali2}, under scale-invariant assumptions on   $\nabla\B$.}. For these reasons, and for the way that we will use (\ref{eq.signass}), 

\begin{center}\emph{we propose that our  assumption (\ref{eq.signass}) generalizes to $L_{\a,V}$ the non-negativity condition on $V$ which is often placed when studying $L_{0,V}$.}
\end{center}

\subsection{Main results for the magnetic Schr\"odinger operator $-(\nabla-i\a)^2+V$}\label{sec.main2a}   Fix a coordinate system for which (\ref{eq.signass}) holds, let
\[
\Sigma_{\m S}{\bf B}(x):=\sum_{(j',k')\in\m S}b_{j'k'},
\]
where $\m S$ is an admissible selection, and let $u=u_{\m S}$ be the  landscape function (with $0$ Dirichlet boundary conditions) from Theorem \ref{thm.exist} for the   Schr\"odinger operator
\begin{equation}\label{eq.electric}
	\tilde L:=-\Delta+\Sigma_{\m S}\B+V,
\end{equation}
so that
\begin{equation}\label{eq.fmu} 
	-\Delta u+\Sigma_{\m S}\B u+Vu=1,\qquad\text{on }\Omega.
\end{equation}

\begin{theorem}[An uncertainty principle for the magnetic Schr\"odinger operator]\label{thm.magup} Let $n\geq2$, and let $\Omega\subset\bb R^n$ be an open connected set, with empty or Lipschitz $(n-1)-$dimensional boundary. Assume that   $\a, {\bf B},V$ verify the assumptions (\ref{eq.bl1}) and (\ref{eq.signass}), and suppose that $\m S$ is an admissible selection such that the Green's function on $\Omega$ for the operator $\tilde L$ from (\ref{eq.electric}) verifies (\ref{eq.greenint}). Then we have that
\begin{equation}\label{eq.up}
\int_\Omega u_{\m S}^2\Big|\nabla\Big(\frac{|f|}{u_{\m S}}\Big)\Big|^2+	\int_\Omega\frac1{u_{\m S}}|f|^2\leq\int_{\Omega}\Big[n|D_{\bf a}f|^2+V|f|^2\Big],
\end{equation}
for any   $f\in C_c^{\infty}(\Omega;\bb C)$.  Here,  $u_{\m S}$ is the  Filoche-Mayboroda landscape function for  $\tilde L$.
\end{theorem}

\begin{remark}\label{rm.dependons} The left-hand side of (\ref{eq.up}) depends on which admissible selection $\m S$ is taken, and in the full generality of the assumption (\ref{eq.signass}),  there may not necessarily exist  a ``maximal'' admissible selection $\m S^*$ such that $u_{\m S^*}$ is pointwise bounded above in $\Omega$ by  $u_{\m S}$, for any other admissible selection  $\m S$. We may remove dependence on $\m S$ from (\ref{eq.up}) by adding up the estimate over all admissible selections $\m S$, of which there are at most a uniform constant depending only on dimension (for a fixed coordinate system). On the other hand, observe that if the stronger direction assumption (\ref{eq.ass}) is verified, then there does exist a maximal admissible selection, due to (\ref{eq.magb2}).
\end{remark}

\begin{remark} The assumption (\ref{eq.signass}) is used to qualitatively assert  the existence and non-negativity of $u$; it might be the case that some further generalization of (\ref{eq.signass}) is possible; for instance, under  subcriticality assumptions on the negative part of $\Sigma_{\m S}\B+V$, but we have not pursued this direction.
\end{remark}
 
The proof of Theorem \ref{thm.magup} is shown in Section \ref{sec.ms}, and relies on the uncertainty principle that $u$ is already known to verify as the Filoche-Mayboroda landscape function for a Schr\"odinger operator with no magnetic field, as well as the identification of the components of the magnetic fields as commutators of the operators (\ref{eq.maggrad}).  This identification or similar ideas have been used by several authors before to obtain uncertainty principles \cite{ahs78, v86, helffer88, shene, ivrii98, nakamura2000,  bls2004,  vt2010,  ctt11}, but usually at the cost of extra derivatives of the magnetic field, or scale-invariant regularity of the direction of the magnetic field. Our   insight is that, under the favored directionality assumption (\ref{eq.signass}), we may use the  a landscape function to dispense of the scale-invariant smoothness condition (\ref{eq.gradb}).  The idea of studying properties of $L_{\a,V}$ by comparison to a non-magnetic Schr\"odinger operator is well-established; indeed, when $V\geq0$ and under minimal integrability conditions on $\a$, the Kato-Simon inequality \cite{kato72, simon79} (see also \cite{ls}) is the semigroup estimate
\[
|e^{-tL_{\a,V}}\varphi|\leq e^{t\Delta}|\varphi|,\qquad\text{for each }\varphi\in C_c^{\infty}(\bb R^n).
\]
In Section \ref{sec.t}, we also show  uncertainty principles  using  landscape functions  in  cases where our directionality assumptions may not apply, or apply only on subdomains of $\bb R^n$.

From Theorem \ref{thm.magup}, a priori exponential decay results can be obtained for  the magnetic Schr\"odinger operator. The  precise definitions of the operators $L_{\a ,V}$, $\n L^{\f a}$ are in Section \ref{sec.ms}. For the a priori exponential decay of the Green's function, we will need to make some more assumptions.

\begin{assumption}\label{ass.g} For $n\geq3$,  assume that\footnote{If $V\geq0$, then $L_{\a,V}$ verifies the properties here stated, as presented in \cite[Section 5]{mp}.} the local weak solutions of $L_{\a,V}\psi=0$ satisfy a scale-invariant Moser estimate\footnote{See Definition 6.1 of \cite{mp}.}, that the Green's function $G_{\a, V}$ for  $L_{\a,V}$ exists on $\Omega$, that it solves $L_{\a,V}G_{\a,V}=0$ outside the diagonal, and that it verifies
\begin{equation}\label{eq.greenbound2}
|G_{\a,V}(x,y)|\lesssim\frac1{|x- y|^{n-2}},\qquad\text{for all }x,y\in\Omega.
\end{equation}
\end{assumption}

\begin{theorem}[Exponential decay of solutions to $L_{\a,V}$]\label{thm.expdecayms} Retain the setting and assumptions of Theorem \ref{thm.magup}, assume that (\ref{eq.nosing}) holds for $\hat u=\hat u_{\m S}$, and moreover, assume that 
\begin{equation}\label{eq.big}
\Sigma_{\m S}\B\geq n V_-=n\max\{0,-V\}.
\end{equation}
Then there exist $\ep>0$ and $C\geq1$, depending only on $n$, so that the following are true.
\begin{enumerate}[(i)]
\item\label{item.lmms}\emph{(Lax-Milgram solutions).} For every  $f\in L^2(\Omega;\bb C)$ with compact support in $\Omega$,
\begin{equation}\label{eq.magup3}\nonumber
\int_\Omega u_{\m S}^2\Big|\nabla\Big(\frac{e^{\ep\rho(\cdot,\supp f,\frac1{\hat u_{\m S}})}|L^{-1}_{\a, V}f|}{u_{\m S}}\Big)\Big|^2+	\int_\Omega\frac1{u_{\m S}}e^{2\ep\rho(\cdot,\supp f,\frac1{\hat u_{\m S}})}|L^{-1}_{\a, V}f|^2\leq C\int_\Omega u_{\m S}|f|^2.
\end{equation} 

\vspace{2mm}
	
\item\label{item.eigenms}\emph{(Eigenfunctions).}  Suppose  there is $\mu>0$ and $\psi\in\n D(\n L^{\f a})$ with $\n L^{\f a}\psi=\mu\psi$. Let
\[
w(x):=\big(\tfrac1{\hat u_{\m S}(x)}-\mu\big)_+=\max\big\{0,\tfrac1{\hat u_{\m S}(x)}-\mu\big\},\quad\text{and}\quad E:=\big\{x\in\Omega: \tfrac1{\hat u_{\m S}(x)}\leq\mu\big\}.
\]
Then  
\begin{equation}\label{eq.eigenms}\nonumber
\int_\Omega u_{\m S}^2\Big|\nabla\Big(\frac{e^{\ep\rho(\cdot,E,w)}|\psi|}{u_{\m S}}\Big)\Big|^2+\int_{\Omega}\Big(\frac1{u_{\m S}}-\mu\Big)_+e^{2\ep\rho(\cdot,E,w)}|\psi|^2  \leq C_n\int_E\Big(\mu-\frac1{u_{\m S}}\Big)_+|\psi|^2.
\end{equation}

\vspace{2mm}

\item\label{item.greenms}\emph{(Green's function).} Suppose that $n\geq3$,  that  $\B_\infty:=\Vert\Sigma_{\m S}\B+V\Vert_{L^{\infty}(\Omega)}<+\infty$, and that Assumption \ref{ass.g} holds. Then there exists $C_G\geq1$, depending on the constants from Assumption \ref{ass.g}, so that for any $x,y\in\Omega$ with
\[
B_x\cap B_y:=B\big(x, \B_\infty^{-1/2}\big)\cap B\big(y,\B_\infty^{-1/2}\big)=\varnothing,\qquad 4B_x\cup4B_y\subset\Omega,
\]
there exist\footnote{See also Remark \ref{rm.neumannms} for a situation where we may use $x,y$ directly in (\ref{eq.greendecayms}).} $\tilde x\in B_x$, $\tilde y\in B_y$ such that 
\begin{equation}\label{eq.greendecayms}
	|G_{\a, V}(\tilde x,\tilde y)|\leq CC_G\B_\infty^{\frac n2}\sqrt{u(\tilde x)}\sqrt{u(\tilde y)}e^{-\ep\rho(\tilde x,\tilde y,\frac1{u_{\m S}})}.
\end{equation}
\end{enumerate}
\end{theorem}

Finally, we restrict to $n\geq3$, and apply Theorem \ref{thm.compare} together with Theorem \ref{thm.magup} to trade the smoothness assumption (\ref{eq.gradb}) for the directionality assumption (\ref{eq.signass}) in the proof of the Fefferman-Phong uncertainty principle for the maximal function $m$.

\begin{corollary}[Uncertainty principle for $L_{\a,V}$ via the Fefferman-Phong-Shen maximal function]\label{cor.upm} Let $n\geq3$, $\Omega=\bb R^n$, and suppose that the operator  $L_{\a,V}$ on $\bb R^n$   verifies the assumptions (\ref{eq.bl1}) and (\ref{eq.signass}). Assume, in addition, that there exists an admissible selection $\m S$ for which $\Sigma_{\m S}\B+V$ is a non-degenerate Shen potential. Then there exists   $C\geq1$, which depends only on $n, C_0, C_1, \delta$, such that for each $f\in C_c^{\infty}(\bb R^n;\bb C)$,
	\begin{equation}\label{eq.upm}
	\int_{\bb R^n} m^2(\cdot,\Sigma_{\m S}\B+V)|f|^2\leq C\int_{\bb R^n}\big[|D_{\bf a}f|^2+V|f|^2\big], 
\end{equation}
where $m(\cdot,\Sigma_{\m S}\B+V)$ is the  function from (\ref{eq.fpsm}).	Moreover, the landscape function $u_{\m S}$ verifies the properties (i), (ii), and (iii) of Corollary \ref{cor.u} (with $\Sigma_{\m S}\B+V$ in place of $V$).
\end{corollary} 

\begin{remark} Note that if $V$ and $\B$ verify (\ref{eq.bl1}), (\ref{eq.ass}), and $|\B|+V$ is a non-degenerate Shen potential, then the assumptions of Corollary \ref{cor.upm} hold, and in this case we may canonically take a maximal admissible selection $\m S$; see Remark \ref{rm.dependons}.
\end{remark}

The proof  is an immediate application of the aforementioned theorems. From Corollary \ref{cor.upm}, we are able to recover several results from the literature without imposing conditions on $\nabla\B$ nor on the pointwise size of $V$, having traded these for the favored directionality assumptions (\ref{eq.signass}) or (\ref{eq.ass}). As a consequence, we resolve Conjecture \ref{conj.s1} in the negative, and we have solved Problem \ref{pro.s2} of extending to higher dimensions the idea of Shen \cite{shenb} to treat the magnetic Schr\"odinger operator with directionality assumptions instead of smoothness assumptions. The  results recovered from the literature by our methods are detailed in  Section \ref{sec.cor}.

 \begin{remark}\label{rm.differentbc} Throughout this paper, we have proved and stated results for the landscape function with homogeneous Dirichlet boundary conditions. However, several of the proofs  in this paper apply also with  natural modifications to problems with different boundary conditions, such as Neumann or mixed boundary conditions, at least in the setting of bounded domains. For the sake of simplicity and brevity, we do not state these results for problems with different boundary conditions.
 \end{remark}

The rest of this paper is organized as follows. In Section \ref{sec.optheory}, we prove Theorem \ref{thm.exist} and give necessary preliminaries. In Section \ref{sec.exp}, we prove Theorems \ref{thm.lmr}, \ref{thm.greenexp}, and \ref{thm.eigen}. In Section \ref{sec.fps}, we prove Theorem \ref{thm.compare}. Finally, in Section \ref{sec.ms}, we prove Theorem \ref{thm.magup}, and give   solutions to Conjecture \ref{conj.s1} and Problem \ref{pro.s2}.

\section{The landscape function on unbounded domains}\label{sec.optheory}

In this section, we provide operator-theoretic preliminaries to construct the landscape function, show  the uncertainty principle, and prove Theorem \ref{thm.exist}. 

\subsection{Preliminaries}\label{sec.prelim} We write $a\lesssim b$ to mean that there exists a uniform constant $C>0$ such that $a\leq Cb$, and $a\approx b$ means that $C^{-1}b\leq a\leq Cb$. Throughout this section, we take the setting and assumptions of Theorem \ref{thm.exist}, unless explicitly mentioned otherwise.

Define the elliptic operator $L$ acting formally on real-valued functions $\psi$ by
\[
L\psi=-\dv(A\nabla\psi)+V\psi=-\sum_{i,j=1}^n\frac{\partial}{\partial x_i}\Big(a_{ij}\frac{\partial\psi}{\partial x_j}\Big)+V\psi.
\]

We will need to consider the weak form of the operator $L$. Recall that $C_c^{\infty}(\Omega)$ is the space of compactly supported smooth functions on $\Omega$,  and   that $W^{1,2}(\Omega)$ is the Sobolev space of square integrable functions on $\Omega$ whose weak derivatives exist in $\Omega$ and are square integrable functions, while $W_0^{1,2}(\Omega)$ is the completion of $C_c^{\infty}(\Omega)$ under the norm
\[
\Vert\psi\Vert_{W^{1,2}(\Omega)}:=\sqrt{\int_{\Omega}|\psi|^2+|\nabla\psi|^2}.
\]
Moreover, the space $\dt W^{1,2}(\Omega)$ consists of the $L^1_{\loc}(\Omega)$ functions whose weak gradient is square integrable over $\Omega$, and for $p\geq1$, we denote by $L^p_c(\Omega)$ the space of (real-valued) Lebesgue $p-$th integrable functions with compact support in $\Omega$.

\begin{definition}[Weak solutions]\label{def.weak} Given  $f\in(C_c^{\infty}(\Omega))^*$, we say that a   function $\psi\in W^{1,2}_{\loc}(\Omega)\cap L^1_{\loc}(\Omega,V\,dx)$  satisfies $L\psi=f$ \emph{in the weak sense} in $\Omega$  if for each $\varphi\in C_c^{\infty}(\Omega)$, we have that
\begin{equation}\label{eq.weakid}
\int_{\Omega}\Big[A\nabla\psi\nabla\varphi+V\psi\varphi\Big]=\langle f,\varphi\rangle.
\end{equation}
\end{definition}

The following local estimate on the gradients of weak solutions is standard.

\begin{lemma}[Caccioppoli inequality]\label{lm.cacc} Retain the setting of Theorem \ref{thm.exist}, let $B$ be a ball such that $2B\subset\Omega$, let $f\in L^q(2B)$ for some $q>\max\{1,n/2\}$, and suppose that $\psi\in W^{1,2}(2B)\cap L^1(2B,V\,dx)$ solves $L\psi=f$ in the weak sense in $2B$. Then there exists $C\geq1$, depending only on $n$ and $\lambda$, such that
	\begin{equation}\nonumber
		\int_B\Big[|\nabla\psi|^2+V\psi^2\Big]\leq C\Big[\frac1{\operatorname{rad}(B)^2}\int_{2B}\psi^2+\int_{2B}f\psi\Big].
	\end{equation}	
\end{lemma}

Define on $C_c^{\infty}(\Omega)$, the (Dirichlet) form
\begin{equation}\label{eq.form}
\f l(\psi,\varphi)=\f l_{\Omega,A,V}(\psi,\varphi)=\int_{\Omega}\Big[A\nabla\psi\nabla\varphi+V\psi\varphi\Big],\qquad\psi,\varphi\in C_c^{\infty}(\Omega).
\end{equation}
Throughout this manuscript, we will omit the subscripts from $\f l$ whenever they are understood from the context. In fact, the function space  on which we have defined the form $\f l$ so far is suboptimal. Owing to the conditions we have placed on $A$ and $V$, the functional
\begin{equation}\label{eq.norm}
\Vert\psi\Vert_{\f l}:=\sqrt{\f l(\psi,\psi)}
\end{equation}
is a norm on $C_c^{\infty}(\Omega)$, and we let $\f D=\f D_{A,V}(\Omega)$ be the completion of $C_c^{\infty}(\Omega)$ under this norm. Thus $\f D$ is complete.

\begin{lemma}[Embeddings of $\f D$ in bounded subdomains]\label{lm.spaces} Retain the setting of Theorem \ref{thm.exist}. The following statements are true.
\begin{enumerate}[(i)]
		\item\label{item.coercive} Let $M$ be any open, bounded subset of $\Omega$ with Lipshitz boundary, with $\int_MV>0$. Then there exists $C$, depending only on $n$, $\lambda$, and $M$, such that
	\begin{equation}\label{eq.coercive}
		\sqrt{\int_{M}\Big[A\nabla\psi\nabla\psi+V\psi^2\Big]}\geq\frac1C \Vert\psi\Vert_{W^{1,2}(M)}
	\end{equation}
	holds for each $\psi\in C^{\infty}(M)$. In particular, if $\psi\in\f D(\Omega)$, then $\psi\in L^p_{\loc}(\Omega)$ for each $p\in[1,\frac{2n}{n-2}]$ if $n\geq3$, and each $p\in[1,\infty)$ if $n=1,2$.
 	\item\label{item.demb} If $\Omega$ is bounded, then  the embedding $\f D(\Omega)\hookrightarrow W_0^{1,2}(\Omega)$ is continuous. Moreover, if $V\in L^{\infty}(\Omega)$, then $\f D(\Omega)= W_0^{1,2}(\Omega)$, and $\Vert\cdot\Vert_{\f l}\approx\Vert\cdot\Vert_{W^{1,2}(\Omega)}$.
\end{enumerate}
\end{lemma}

\noindent\emph{Proof.} We prove \ref{item.coercive} first.  By the ellipticity of $A$ and the definition of the norm $\Vert\cdot\Vert_{\f l}$, we have that
\[
\Vert\nabla\psi\Vert_{L^2(M)}^2=\int_{M}|\nabla\psi|^2\leq\lambda\int_{M}A\nabla\psi\nabla\psi 
\]
for each $\psi\in C^{\infty}(\Omega)$. Then, the fact that $\Vert\psi\Vert_{L^2(M)}$ is controlled by the left-hand side of (\ref{eq.coercive}) can be proved via a standard compactness argument using the Rellich-Kondrachov theorem\footnote{When $n=1$, we may not use this theorem, but instead we appeal to Morrey's inequality and the Arzelà-Ascoli theorem to reach the same conclusion.}, where we must apply the fact that $V$ is positive on a subset of positive measure. As for \ref{item.demb},   the trace $\psi|_{\partial\Omega}$ vanishes via the trace theorem for $C_c^{\infty}(\Omega)$ functions and a completion argument, and therefore $\psi\in L^2(\Omega)$ due to the Poincar\'e-Sobolev inequality when $\Omega$ is bounded. We omit further details. \hfill{$\square$}

Note that the form $\f l$ is bounded and coercive on $\f D$ (the coercivity follows by definition and the non-degeneracy of $V$, while the boundedness follows from the Cauchy-Schwartz inequality and (\ref{eq.elliptic})), which turns $\f D$ into a Hilbert space. Denote by $\f D^*$ the space of bounded linear functionals on $\f D$, and note that for each $q>\max\{1,\frac{2n}{n+2}\}$, the space $L^q_c(\Omega)$  is contained in $\f D^*$ by virtue of Lemma \ref{lm.spaces} \ref{item.coercive} and the Sobolev embedding theorems. Define the \emph{Dirichlet operator} $L=L_{\Omega,A,V}:\f D\ra\f D^*$ in the following way: for each $\psi\in\f D$, $L\psi\in\f D^*$ is the functional given by
\[
\langle L\psi,\varphi\rangle:=\f l(\psi,\varphi),\qquad\text{for each }\varphi\in\f D.
\]
We may thus apply the Lax-Milgram theorem to see that $L:\f D\ra\f D^*$ is bounded and invertible. Consequently, for each $f\in\f D^*$, there exists a unique $\psi=L^{-1}f\in\f D$ such that
\[
\f l(\psi,\varphi)=\langle f,\varphi\rangle,\qquad\text{for each }\varphi\in\f D.
\]

We now state the  weak maximum principle, proven in \cite[Proposition 3.2]{adfjm} for the case of $V\in L^{\infty}(\Omega)$ and symmetric $A$; we omit the proof of the general case as it is very similar. If $f\in\f D^*$, we say that $f\geq0$ if $\langle f,\varphi\rangle\geq0$ for each $\varphi\geq0$ in $\f D$.

\begin{lemma}[Weak maximum principle, \cite{adfjm}]\label{lm.maxprin} Retain the setting of Theorem \ref{thm.exist}.    Suppose that $\psi$ is a continuous function such that $\Vert\psi\Vert_{\f l}<\infty$, that $\psi\geq0$ on $\partial\Omega$, and that  $L\psi=f\geq0$ in the weak sense in $\Omega$, for $f\in\f D^*$. Then $\psi\geq0$ in $\Omega$.
\end{lemma}

We will repeatedly make use of the following   technical fact about approximations of solutions; its proof is very similar to that of \cite[Lemma 5.31]{mp}, and thus omitted.

\begin{lemma}[Approximation of solutions via bounded potentials]\label{lm.approx}  Retain the setting of Theorem \ref{thm.exist}. Fix $f\in L^q_c(\overline{\Omega})$  for some $q>\max\{1,n/2\}$, and for each $N\in\bb N$, let $V_N:=\min\{V,N\}$, $\psi_n:=L_{V_N}^{-1}f$, $\psi:=L_V^{-1}f$. Then $\nabla\psi_N\ra\nabla\psi$ strongly in $L^2(\Omega)^n$, $V_N^{1/2}\psi_N\ra V^{1/2}\psi$ strongly in $L^2(\Omega)$, and $\psi_N\ra\psi$  strongly in $L^2_{\loc}(\overline{\Omega})$.
	
Moreover, if $f\geq0$, then $\psi_N\searrow\psi$ pointwise a.e. in $\Omega$.
\end{lemma}

The last lemma allows us to upgrade our weak maximum principle.

\begin{corollary} \label{cor.maxprin1} Retain the setting of Theorem \ref{thm.exist}, and let $0\leq f\in L^q_c(\overline{\Omega})$ for some $q>\max\{1,n/2\}$. Then $L^{-1}f\geq0$ in $\Omega$. 
\end{corollary}

When the singularities of $V$ are controlled, the following stronger version of the maximum principle is classical.

\begin{lemma}[Strong maximum principle]\label{prop.strong} Retain the setting of Theorem \ref{thm.exist},    assume that $V\in L^{q}(\Omega)$ for some $q>\max\{1,n/2\}$, and that $\Omega$ is bounded. Then, if $\psi\in\f D$ verifies $L\psi=f\geq0$, $f\in L^{\infty}(\Omega)$, with $\int_\Omega f>0$, it follows that $\psi>0$ in $\Omega$.
\end{lemma}

The weak maximum principle can be used to verify the existence of Green's function $G$ as an integral kernel of the operator $L^{-1}$. We state the definition of $G$ below, and defer its existence to Proposition \ref{prop.green} in Appendix \ref{sec.green}. 

\begin{definition}[Green's function]\label{def.green} We say that a measurable function
\[
G:\Omega\times\Omega\ra\bb R
\]
is the \emph{Green's function} for the operator $L$ on $\Omega$ if the following statements hold.
\begin{enumerate}[(i)]
	\item\label{item.greenloc} For each $x\in\Omega$, $0\leq G(x,\cdot)\in L^p_{\loc}(\Omega)$ for some $p>1$.
	\item\label{item.greenres} For each $f\in L^{\infty}_c(\Omega)$, we have the identity
	\begin{equation}\label{eq.green}
	(L^{-1}f)(x)=\int_{\Omega}G(x,y)f(y)\,dy\qquad\text{for a.e. }x\in\Omega.
	\end{equation}  
\end{enumerate}
\end{definition}

When $\Omega=\bb R^n$, the Green's function is known as the \emph{fundamental solution}. Note that, under   stronger assumptions on $V$ (say, for instance, $0\leq V\in L^{\frac n2+\ep}_{\loc}(\Omega)$ when $n\geq3$), it can be shown \cite{dhm} that   for each $\varphi\in C_c^{\infty}(\overline{\Omega})$, the identity
\begin{equation}\label{eq.greensol}
	\int_{\Omega}\Big[A(x)\nabla G(x,y)\nabla\varphi(x)+V(x)G(x,y)\varphi(x)\Big]\,dx=\varphi(y)
\end{equation}
holds for all $y\in\Omega$.  However, for our singular potentials,   (\ref{eq.greensol}) may not be true,  as remarked in Section \ref{sec.main}.

\subsection{Construction of the landscape function and the uncertainty principle}

\begin{definition}[The landscape function on bounded domains]\label{def.landscape} If $\Omega$ is bounded, let $u=u_{\Omega}:=L_{\Omega}^{-1}\1_{\Omega}$. Thus,  $u$ is the unique weak solution to the Dirichlet problem $Lu=1$ with Dirichlet data $0$, called the \emph{landscape function}. Observe that $\Tr u|_{\partial\Omega}\equiv0$.
\end{definition}

The following lemma is crucial in understanding the usefulness of the landscape function. It is shown in \cite{adfjm} for symmetric matrices, under the additional assumption that $f\in L^{\infty}(\Omega)$, and for a landscape function which is strongly positive on $\Omega$. Our setting is a bit more complicated; let us give a careful proof.

\begin{lemma}[Operator conjugation by $\frac1u$ for bounded potentials]\label{lm.conj} Retain the setting of Theorem \ref{thm.exist}, and moreover, assume that $V\in L^{\infty}(\Omega)$, and that $\Omega$ is bounded. Let $u=L^{-1}_{\Omega}\1_{\Omega}$. Then for each $f\in W_0^{1,2}(\Omega)$, the identity
	\begin{multline}\label{eq.conj}
		\int_\Omega\Big[A\nabla f\nabla f+Vf^2\Big]=\int_\Omega\Big[u^2A\nabla(f/u)\nabla(f/u)+\frac1uf^2\Big]\\+\int_\Omega\frac fu\Big[A\nabla f\nabla u-A\nabla u\nabla f\Big] 
	\end{multline}
holds. In particular, if $A$ is symmetric, the last integral vanishes for any $f$.
\end{lemma}

\begin{remark} Since  $\frac1u$ necessarily blows up near $\partial\Omega$,  the finiteness of the integrals in the right-hand side of (\ref{eq.conj}) is not trivial.
\end{remark}

\noindent\emph{Proof of Lemma \ref{lm.conj}.} {\bf Step 1: Case of $f\in C_c^{\infty}(\Omega)$.}   Since $V\in L^{\infty}(\Omega)$ and $\Omega$ is bounded, then $\f D=W_0^{1,2}(\Omega)$ by Lemma \ref{lm.spaces}, and by the H\"older continuity of $u$ and the strong maximum principle, we have that $\frac1u\in L^{\infty}(\supp f)$. Then it is easy to see that $f^2/u$ belongs to $\f D$ since $\frac1u\in L^{\infty}(\supp f)$ and $\nabla u\in L^2(\Omega)$. Thus we may use $f^2/u$ as a test function in the identity $\f l(u,\varphi)=\langle1,\varphi\rangle$ to see that
\begin{equation}\label{eq.testu}
\int_{\Omega}\Big[A\nabla u\nabla(f^2/u)+Vf^2\Big]=\int_{\Omega}\frac1uf^2.
\end{equation}
Using the product rule we may rewrite the first term of (\ref{eq.testu}) as follows:
\begin{equation}\label{eq.rewrite}
A\nabla u\nabla(f^2/u)=A\nabla f\nabla f-u^2A\nabla(f/u)\nabla(f/u)+\frac fu\Big[A\nabla u\nabla f-A\nabla f\nabla u\Big].
\end{equation}
Putting the last two identities together yields (\ref{eq.conj}) in this case.

{\bf Step 2: Finiteness of the integrals.} We now show that (\ref{eq.conj}) holds for $f\in W_0^{1,2}(\Omega)$. Let $\{f_k\}_{k=1}^{\infty}$ be a family of functions in  $C_c^{\infty}(\Omega)$ such that $f_k\ra f$ strongly in $W_0^{1,2}(\Omega)$. In particular, there is a subsequence where $f_{k'}\ra f$, $\nabla f_{k'}\ra\nabla f$ pointwise a.e. in $\Omega$ as $k'\ra\infty$, and for simplicity of notation we now take this subsequence to be the whole sequence. Observe that (\ref{eq.conj}) is true with $f$ replaced by $f_k$, for each $k\in\bb N$.  By using the Cauchy inequality with $\ep>0$ and (\ref{eq.elliptic}), we easily obtain the pointwise estimate
\[
\Big|2\frac{f_k}uA\nabla u\nabla f_k\Big|\leq\frac{f_k^2}{u^2}A\nabla u\nabla u+\frac{1}{\lambda^4}A\nabla f_k\nabla f_k.
\]
Since 
\begin{equation}\label{eq.pr}
	A\nabla u\nabla(f_k^2/u)=2\frac{f_k}uA\nabla u\nabla f_k-\frac{f_k^2}{u^2}A\nabla u\nabla u,
\end{equation}
it therefore follows that 
\[
A\nabla u\nabla(f_k^2/u)\leq\frac{1}{\lambda^4}A\nabla f_k\nabla f_k.
\]
We now plug this last estimate in (\ref{eq.testu}) to see that
\begin{equation}\label{eq.up2}
\int_{\Omega}\frac1uf_k^2\leq\int_{\Omega}\Big[\frac1{\lambda^4}A\nabla f_k\nabla f_k+Vf_k^2\Big],\qquad\text{for each }k\in\bb N.
\end{equation}
Let us point out that for any $j,\ell\in\bb N$, $f_j-f_{\ell}\in C_c^{\infty}(\Omega)$, and thus we may obtain (\ref{eq.up2}) with $f_k$ replaced by $f_j-f_{\ell}$. Since $\{f_k\}$ is Cauchy in $W_0^{1,2}(\Omega)$, this implies that $\{f_k\}$ is Cauchy in $L^2(\Omega,\frac1u\,dx)$. It follows that $\frac1uf^2\in L^1(\Omega)$, and $\frac1uf_k^2\ra\frac1uf^2$ in $L^1(\Omega)$. 

Now we investigate the finiteness of the rest of the integrals in (\ref{eq.conj}).  For each $k\in\bb N$, we use (\ref{eq.testu}) and (\ref{eq.pr}) to  obtain the identity
\begin{equation}\label{eq.conjk}
\int_{\Omega}\Big[\frac{f_k^2}{u^2}A\nabla u\nabla u-2\frac{f_k}uA\nabla u\nabla f_k\Big]=\int_\Omega\Big(V-\frac1u\Big)f_k^2.
\end{equation}
The right-hand side of (\ref{eq.conjk}) has a limit as $k\ra\infty$, whence it must be the case that the left-hand side also has a limit. In particular, there exists $C\in[0,\infty)$ such that
\[
\int_{\Omega}\Big[\frac{f_k^2}{u^2}A\nabla u\nabla u-2\frac{f_k}uA\nabla u\nabla f_k\Big]\leq C,\qquad\text{for all }k\in\bb N.
\]
Now we may use the Cauchy inequality with $\ep>0$ to see that
\[
\Big|\int_{\Omega}2\frac{f_k}uA\nabla u\nabla f_k\Big|\leq\frac12\int_{\Omega}\frac{f_k^2}{u^2}A\nabla u\nabla u~+~\frac{2}{\lambda^3}\int_{\Omega}|\nabla f_k|^2,\qquad\text{for all }k\in\bb N.
\]
Therefore the estimate
\begin{equation}\label{eq.calc1}
\int_{\Omega}\frac{f_k^2}{u^2}A\nabla u\nabla u\leq2C+\frac{4}{\lambda^3}\int_{\Omega}|\nabla f_k|^2\leq C'
\end{equation}
holds for all $k\in\bb N$, and $C'>0$ is uniform in $k$. By Fatou's Lemma it follows that $\frac{f^2}{u^2}A\nabla u\nabla u\in L^1(\Omega)$, and similarly we obtain that $\frac fuA\nabla u\nabla f, \frac fuA\nabla f\nabla u\in L^1(\Omega)$. Then, by (\ref{eq.rewrite}) and (\ref{eq.pr}), we also conclude that $u^2A\nabla(f/u)\nabla(f/u)\in L^1(\Omega)$, and thus every integral in (\ref{eq.conj}) is finite. 

{\bf Step 3: Passing to the limit.} It remains to pass the identity (\ref{eq.conj}), currently shown for $f$ replaced by $f_k$, to the limit as $k\ra\infty$.  To this end, we first show that $\frac{f_k}u|\nabla u|\nabla f_k|\ra\frac fu|\nabla u||\nabla f|$ in $L^1(\Omega)$. Since $\frac fu\nabla u\in L^2(\Omega)$, $\{\frac{f_k}u|\nabla u|\}$ is uniformly bounded in $L^2(\Omega)$, and
\begin{multline}\nonumber
	\Big|\int_{\Omega}\frac{|f_k|}u|\nabla u||\nabla f_k|-\int_{\Omega}\frac{|f|}u|\nabla u||\nabla f|\Big|\leq\int_{\Omega}\frac{|f_k|}u|\nabla u||\nabla f_k-\nabla f|~+~\Big|\int_{\Omega}\frac{|f_k|-|f|}u|\nabla u||\nabla f|\Big|\\  \leq\Big\Vert\frac{|f_k|}u|\nabla u|\Big\Vert_{L^2(\Omega)}\Vert\nabla f_k-\nabla f\Vert_{L^2(\Omega)}~+~\Big|\int_{\Omega}\frac{|f_k|-|f|}u|\nabla u||\nabla f|\Big|,
\end{multline}
by \cite[Exercise 2.21]{folland}, we reach the desired conclusion if we prove that $\frac{|f_k|}u|\nabla u|\ra\frac{|f|}u|\nabla u|$ weakly in $L^2(\Omega)$. But if $\varphi\in C_c^{\infty}(\Omega)$, then $\frac1u|\nabla u|\in L^2(\supp\varphi)$ and therefore
\[
\Big|\int_{\Omega}\frac{|f_k|-|f|}u|\nabla u|\varphi\Big|\leq\Vert f_k-f\Vert_{L^2(\Omega)}\Vert\frac1u|\nabla u|\varphi\Vert_{L^2(\supp\varphi)}\longrightarrow0
\]
as $k\ra\infty$, whence $\frac{|f_k|}u|\nabla u\ra\frac{|f|}u|\nabla u|$ weakly in $L^2(\Omega)$, and   our claim is shown. Then, by the Lebesgue Dominated Convergence Theorem and the boundedness of $A$, we have that
\[
\frac{f_k}uA\nabla u\nabla f_k\ra\frac fuA\nabla u\nabla f,\qquad\frac{f_k}uA\nabla f_k\nabla u\ra\frac fuA\nabla f\nabla u
\]
in $L^1(\Omega)$ as $k\ra\infty$. Now we need only pass $I(f_k):=u^2A\nabla(f_k/u)\nabla(f_k/u)$ to the limit in $L^1(\Omega)$. But in fact, by using $f_k-f_\ell$ in place of $f_k$ in (\ref{eq.conj}), we deduce that $0\leq I(f_k-f_\ell)=\n F(f_k-f_\ell)$, where $\n F(f_k-f_\ell)\ra0$ as $k,\ell\ra\infty$, since we have already shown every other term in (\ref{eq.conj}) can pass to the limit in $L^1(\Omega)$. The desired result follows.\hfill{$\square$}

The previous lemma allows us to show the uncertainty principle that the landscape function satisfies, for singular potentials.

\begin{corollary}[Uncertainty Principle on bounded domains]\label{cor.uncertainty}  Retain the setting of Theorem \ref{thm.exist}, and assume that $\Omega$ is bounded. Then $\frac1u\in L^1_{\loc}(\Omega)$, $\nabla\log u\in L^2_{\loc}(\Omega)$, and for each $f\in\f D(\Omega)$, each of the integrands in (\ref{eq.conj}) lies in $L^1(\Omega)$, and we have that
\begin{equation}\label{eq.unprinciple}
\int_{\Omega}\frac1u f^2\leq\int_{\Omega}\frac1{\lambda^4}A\nabla f\nabla f+\int_{\Omega}Vf^2.
\end{equation} 
If we further assume that $A$ is symmetric, then   for each $f\in\f D(\Omega)$, we have that
\begin{equation}\label{eq.unprinciplesym}
	\int_{\Omega}\frac1u f^2+\int_\Omega u^2A\nabla(f/u)\nabla(f/u)\leq\int_{\Omega}A\nabla f\nabla f+\int_{\Omega}Vf^2.
\end{equation}  
\end{corollary}

\noindent\emph{Proof.} If $V_N:=\min\{V,N\}$, then $0\leq u\swarrow u_N:=L_{V_N}^{-1}\1_{\Omega}$ pointwise a.e. in $\Omega$, by virtue of Lemma \ref{lm.approx} and Corollary \ref{cor.maxprin1}. Then the estimate (\ref{eq.unprinciple}) follows by Fatou's lemma and the estimate (\ref{eq.up2}) shown in the proof of Lemma \ref{lm.conj}. Since (\ref{eq.unprinciple}) holds for arbitrary $f\in\f D(\Omega)\subset W_0^{1,2}(\Omega)$, we see that $\frac1u\in L^1_{\loc}(\Omega)$. Moreover, from the proof of Lemma \ref{lm.conj}, we may similarly check that (\ref{eq.calc1}) implies $\nabla\log u\in L^2_{\loc}(\Omega)$, and from there, the finiteness of each of the integrals in (\ref{eq.conj}) follows. Finally, if $A$ is symmetric, then (\ref{eq.unprinciplesym}) follows directly from (\ref{eq.conj}) and Fatou's lemma.\hfill{$\square$}

We turn to the proof of our first main result.  

\noindent\emph{Proof of Theorem \ref{thm.exist}.} The statement \ref{item.u} follows from the weak maximum principle. 

Proof of \ref{item.equiv} ($\impliedby$). Suppose   that (\ref{eq.greenint}) holds. Fix $\varphi\in C_c^{\infty}(\Omega)$, let $x_0\in\Omega$ be as in \ref{item.u}, and let $R>R_0$ with $R_0$   large enough so that $\supp\varphi\subset\Omega_{\lfloor R_0/2\rfloor}$. By the weak maximum principle and the Green's representation formula (\ref{eq.green}), it is easy to see that for a.e. $x\in\Omega_R$, $G_R(x,y)\leq G(x,y)$ for a.e. $y\in\Omega_R$. Let $B$ be any ball such that $2B\subset\Omega_{R_0}$. By the De Giorgi-Nash-Moser estimate\footnote{If $n=1$, we get the same conclusion by Morrey's inequality.} \cite[Theorem 4.1]{hanlin}, we have that for all $R>R_0$,
\begin{multline}\label{eq.essbdd}
\Vert u_R\Vert_{L^{\infty}(B)}\leq C_B\Big[\Vert u_R\Vert_{L^q(2B)}+1\Big]=C_B\Big[\Big\Vert\int_{\Omega_R}G_R(x,y)\,dy\Big\Vert_{L^q(2B,\,dx)}+1\Big]\\ \leq C_B\Big[ \Big\Vert\int_{\Omega}G(x,y)\,dy\Big\Vert_{L^q(2B,\,dx)}+1\Big]\leq\tilde C_B
\end{multline}
where we have used (\ref{eq.greenint}), and the constant $\tilde C_B$ does not depend on $R$. Since $\supp\varphi$ is compact, we may cover it by finitely many balls $B$ of bounded radius $\rad(B)\approx1$ and such that    $2B\subset\Omega_{R_0}$. Then, from (\ref{eq.essbdd}) we may furnish  the estimate
\begin{equation}\label{eq.bdd}
\Vert u_R\Vert_{L^{\infty}(\supp\varphi)}\leq C_{\varphi},\qquad\text{for each }R>R_0.
\end{equation}
This estimate   implies that $u$, defined via (\ref{eq.limit}),   is essentially bounded in $\supp\varphi$. With (\ref{eq.bdd}) at hand and via a similar covering argument, the Caccioppoli inequality Lemma \ref{lm.cacc} gives that $\{\nabla u_R\}_{R>R_0}$ is a uniformly bounded sequence in $L^2(\supp\varphi)$, so that a subsequence   converges weakly in $L^2(\supp\varphi)$, say, to ${\bf g}\in L^2(\supp\varphi)^n$. Then we must have $\nabla u={\bf g}$, and thus we must have that the whole sequence $\nabla u_R$ converges weakly to $\nabla u$ in $L^2(\supp\varphi)$. Moreover, by (\ref{eq.limit}), (\ref{eq.bdd}), and the LDCT, we have that $Vu_R\varphi\ra Vu\varphi$ in $L^1(\supp\varphi)$. Since $u_R$ solves $Lu_R=1$ weakly in $\Omega_{R_0}$ for each $R>R_0$, we may use the test function $\varphi$ for all $R>R_0$ and pass to the limit $R\ra\infty$ in the integral identity to see that
\[
\int_{\Omega}\Big[A\nabla u\nabla\varphi+Vu\varphi\Big]=\int_\Omega\varphi.
\]
But because $\varphi\in C_c^{\infty}(\Omega)$ was arbitrary, we have obtained the desired conclusion.

Proof of \ref{item.equiv} ($\implies$). Assume that the function $u$ of (\ref{eq.limit}) lies in $W^{1,2}_{\loc}(\Omega)\cap L^2_{\loc}(\Omega,V\,dx)$ and solves the equation $Lu=1$ in the weak sense in $\Omega$. Extend each $G_R$ by $0$ outside of its original domain of definition $\Omega_R\times\Omega_R$, and observe that by the weak maximum principle, for a.e. $x\in\Omega$, $\{G_R(x,y)\}_{R>1}$ is a monotone non-decreasing sequence of functions in $y$. Hence, for a.e. $x\in\Omega$, let $J(x,y):=\lim_{R\ra\infty}G_R(x,y)=\limsup_{R\ra\infty}G_R(x,y)$, and by the weak maximum principle again it is easy to see that $0\leq J(x,y)\leq G(x,y)$.

Fix a ball $B$ such that $2B\subset\Omega$, and let $R_0$ be so large that $2B\subset\Omega_{R_0}$. Since $Lu=1$ in the weak sense in $\Omega$, then by the Moser estimate we have that $\Vert u\Vert_{L^{\infty}(B)}\leq C_B$. The Lebesgue Monotone Convergence Theorem lets us conclude that 
\[
\int_{\Omega}J(x,y)\,dy=\lim_{R\ra\infty}\int_{\Omega_R}G_R(x,y)\,dy=u(x),\qquad\text{for a.e. }x\in B,
\]
whence $\int_{\Omega}J(x,y)\,dy\in L^{\infty}(B,\,dx)$. It remains only to check that $J\equiv G$ as a measurable function on $\Omega\times\Omega$. But, by considering arbitrary $f\in L^{\infty}_c(\Omega)$, it may be checked, similarly as in the convergence argument in the proof of \ref{item.equiv} ($\impliedby$) (see also Lemma \ref{lm.approx} and \cite[Lemma 5.31]{mp}), that $\int_{\Omega}J(x,y)f(y)\,dy=L^{-1}_{\Omega}f$. From this fact and the uniqueness of the Green's function, we conclude that $J\equiv G$. Then (\ref{eq.greenint}) holds with $q=\infty$.

Proof of \ref{item.under}.  Suppose that (\ref{eq.greenint}) holds. In the proof of \ref{item.equiv} ($\implies$) we already saw that $u=\int_\Omega G(\cdot,y)\,dy$ a.e. on $\Omega$. If $\Omega$ is bounded, by the uniqueness of the Green's function there is no conflict with Definition \ref{def.landscape} of the landscape function on bounded domains. The uncertainty principles follow easily by Corollary \ref{cor.uncertainty} (and its proof), Fatou's lemma and the convergence properties of $u_R$ to $u$. We omit further details.\hfill{$\square$}

We finish this section with a lower bound property for the landscape function in terms of the original potential $V$. When $V$ and $\Omega$ are bounded, for a Neumann landscape function $u_{\f N}$ it is easy to obtain the estimate $u_{\f N}\geq\Vert V\Vert_{L^{\infty}(\Omega)}^{-1}$  via the maximum principle \cite[Proposition 3.2]{adfjm}, but since the Dirichlet landscape function vanishes at the boundary, such an estimate cannot hold over $\Omega$ for our $u$. Nevertheless, it is still possible to find plenty of balls where the estimate holds up to a uniform multiplicative constant. We will use the following result heavily in the proof of the exponential decay estimate for Green's function for bounded potentials.

\begin{proposition}\label{prop.harnacku} Retain the setting of Theorem \ref{thm.exist},  assume that $n\geq2$, that (\ref{eq.greenint}) holds, and that $V\in L^{\infty}(\Omega)$.  Then there exist   constants $c_1,c_2\in(0,1)$, depending only on $n$ and $\lambda$, such that for any $x_0\in\Omega$ with $B(x_0,4\Vert V\Vert_{L^{\infty}(\Omega)}^{-1/2})\subset\Omega$, there exists $\tilde x\in B_0:=B(x_0,\Vert V\Vert_{L^{\infty}(\Omega)}^{-1/2})$ so that
\begin{equation}\label{eq.lbound}
u(x)\geq c_1\Vert V\Vert_{L^{\infty}(\Omega)}^{-1} \quad\text{for every }x\in B(\tilde x,c_2\Vert V\Vert_{L^{\infty}(\Omega)}^{-1/2}).
\end{equation}
\end{proposition}

\noindent\emph{Proof.} We first show that there exists a point $\tilde x\in B_0$ so that $u(\tilde x)>c\Vert V\Vert_{L^{\infty}(\Omega)}^{-1}$, for $c$ small enough. Fix a parameter $\mu\in(0,1)$, and suppose that $u\leq\mu\Vert V\Vert_{L^{\infty}(\Omega)}^{-1}$ on $B_0$.  Fix any $f\in C_c^{\infty}(B_0)$, and by the uncertainty principle (\ref{eq.unprinciple}) we obtain that
\begin{equation}\label{eq.valid}
\int_{B_0}|\nabla f|^2\geq\lambda^5\Big[\frac1{\mu}-1\Big]\Vert V\Vert_{L^{\infty}(\Omega)}\int_{B_0}|f|^2 .
\end{equation}
In particular, (\ref{eq.valid}) must be valid for the cut-off function $f$ given as $f\equiv1$ on $\frac12B_0$, $0\leq f\leq1$, and $|\nabla f|\leq2\Vert V\Vert_{L^{\infty}(\Omega)}^{1/2}$. From (\ref{eq.valid}), it follows that $\mu^{-1}\leq 1+2^{n+1}\lambda^{-5}$. Hence, if $\mu=2^{-n-2}\lambda^5$, then there must exist $\tilde x\in B_0$ so that $u(\tilde x)>\mu\Vert V\Vert_{L^{\infty}(\Omega)}^{-1}$.

To finish the proof, we fix $\delta\in(0,1)$ and recall the Harnack inequality \cite[Theorems 8.17 and 8.18]{gtpde} that $u$ satisfies:
\begin{equation}\label{eq.harnack1}
\sup_{x\in B(\tilde x, \delta\Vert V\Vert_{L^{\infty}(\Omega)}^{-1/2})}u(x)\leq\tilde C\Big[\inf_{x\in B(\tilde x, \delta\Vert V\Vert_{L^{\infty}(\Omega)}^{-1/2})}u(x)+\frac{\delta^2}{\Vert V\Vert_{L^{\infty}(\Omega)}}\Big],
\end{equation}
where $\tilde C$ depends only on $n$ and $\lambda$, since the dependence on $(\Vert V\Vert_{L^{\infty}(\Omega)}^{1/2})(\Vert V\Vert_{L^{\infty}(\Omega)}^{-1/2})=1$ disappears\footnote{Technically, the estimates as written in \cite{gtpde} allow a dependence on $\delta$; but as long as $\delta\in(0,1)$, this dependence can be ignored because the constant only grows when $\delta>1$; see also the discussion after Theorem 8.20 in \cite{gtpde}.}. It follows that if $\delta=\sqrt{\mu/2\tilde C}$, then for any $x\in B(\tilde x, \delta\Vert V\Vert_{L^{\infty}(\Omega)}^{-1/2})$,
\[
u(x)\geq\frac{\mu}{\tilde C}\frac1{\Vert V\Vert_{L^{\infty}(\Omega)}}-\frac{\delta^2}{\Vert V\Vert_{L^{\infty}(\Omega)}}\geq\frac{\mu}{2\tilde C}\frac1{\Vert V\Vert_{L^{\infty}(\Omega)}},
\]
as desired. \hfill{$\square$}

\section{A priori exponential decay results for $-\dv A\nabla+V$}\label{sec.exp}

Throughout this section, we take the setting of Theorem \ref{thm.exist}, and additionally, we assume that (\ref{eq.greenint}) holds, so that the landscape function $u$, as defined in (\ref{eq.limit}), makes sense over $\Omega$ as a solution to the equation $Lu=1$.

Let $w$ be a non-negative, continuous function on $\Omega$. Given an elliptic matrix $\m A$, we denote the entries of $\m A^{-1}$ by $\m A^{-1}_{ij}(x)$. We define the distance $\rho_{\m A}(x,y,w)$ on $\Omega$ for the degenerate Riemannian metric $ds^2=w(x)\sum\m A^{-1}_{ij}dx_idx_j$ by
\begin{equation}\label{eq.agmon}
	\rho_{\m A}(x,y,w)=\inf_{\gamma}\int_0^1\Big(w(\gamma(t))\sum_{i,j=1}^n\m A^{-1}_{ij}(\gamma(t))\dot\gamma_i(t)\dot\gamma_j(t)\Big)^{1/2}\,dt,
\end{equation}
if $\m A$ is symmetric and continuous, where the infimum is taken over all absolutely continuous paths $\gamma:[0,1]\ra\Omega$ such that $\gamma(0)=x$ and $\gamma(1)=y$. If the matrix $\m A$ is not symmetric or not continuous, then we will use the distance $\rho(x,y,w):=\rho_I(x,y,w)$, where $I$ is the $n\times n$ identity matrix. The main reason to take the assumption of continuity when $\m A$ is symmetric is that it will allow us a quantitative bound with no dependence on the ellipticity constants. If $\m A$ is not symmetric, we have already introduced quantitative dependence on the ellipticity constants in Corollary \ref{cor.uncertainty}, so in this case we do not bother with the continuity assumption.

\begin{lemma}[\cite{agmon}, Theorem 4, p.18]\label{lm.lip} Suppose that $\m A$ is a, symmetric, continuous, elliptic matrix, and let $w$ be a non-negative, continuous function on $\Omega$. If $\phi$ is real-valued and $|\phi(x)-\phi(y)|\leq\rho_{\m A}(x,y,w)$ for all $x,y\in M$, then  $\phi$ is locally Lipschitz, and 
\begin{equation}\label{eq.phiest}
\m A(x)\nabla\phi(x)\nabla\phi(x)\leq w(x),\qquad\text{for each }x\in M.
\end{equation} 
In particular, we may choose $\phi=\rho(\cdot,E,w):=\inf_{y\in E}\rho(\cdot,y,w)$ for any nonempty  $E\subset\Omega$. 
\end{lemma}

\begin{remark}\label{rm.which} If $V\in L^{\infty}(\Omega)$, then the landscape function $u$ is locally H\"older continuous, and positive on $\Omega$ by virtue of the strong maximum principle, and so in the case of bounded potentials, we may take $w=1/u$ with no issues. However, if we merely have $V\in L^1_{\loc}(\Omega)$, then we may lose both the positivity and continuity of $u$. In fact, in this general case, we may have different representatives of $u$ which agree pointwise a.e. on $\Omega$, but their respective Agmon distances may be quite different. To fix this issue, we will always consider the representative $\hat u$ pointwise given by
\[
\hat u(x):=\lim_{N\ra\infty}u_N(x),\qquad\text{for each }x\in\Omega,
\]
where for each $N\in\bb N$, $u_N$ is the (Dirichlet) landscape function on $\Omega$ for the operator $-\dv A\nabla+V_N$, $V_N:=\min\{V,N\}$, as constructed in Theorem \ref{thm.exist}. The limit exists pointwise  because $\{u_N\}$ is a monotone non-increasing family of functions on $\Omega$. That $\hat u$ solves $L\hat u=1$ on $\Omega$ and agrees with $u$ as an element of $W^{1,2}_{\loc}(\Omega)$ (where $u$ is as defined in Theorem \ref{thm.exist}) is a consequence of our approximation results from Section \ref{sec.optheory}. With this definition, we will let $w=1/\hat u$. Of course,   if $V\in L^{\infty}(\Omega)$, then $u\equiv\hat u$ pointwise on $\Omega$.
\end{remark}

\noindent\emph{Proof of Theorem \ref{thm.lmr}}. Fix $f\in L^2_c(\Omega)$, and let us first show the estimate (\ref{eq.lmest}) when $V\in L^{\infty}(\Omega)$. Write $\psi:=L^{-1}f$, and for each $M\in\bb N$, let $\phi_M:=\min\{M,\rho(x,\supp f,\frac1u)\}$ and $g=e^{\ep\phi_M}$ for some $\ep>0$. Then $g\psi\in\f D(\Omega)$: indeed, $g\psi\in L^2(\Omega,V\,dx)$ because $g\leq M$ on $\Omega$, and since $\nabla(g\psi)=\ep\psi e^{\ep\phi_M}\nabla\phi+e^{\ep\phi_M}\nabla\psi$, then $\nabla(g\psi)\in L^2(\Omega)$, because $\nabla\psi\in L^2(\Omega)$  and $|\nabla\phi|\psi\lesssim\frac1{\sqrt u}|\psi|\in L^2(\Omega)$ by (\ref{eq.unprinciple}). Using (\ref{eq.unprinciple}) with $g\psi$, we see that
\begin{equation}\label{eq.exp1}
\int_\Omega\frac1ug^2|\psi|^2\leq C\int_\Omega\Big[ A\nabla(g\psi)\nabla(g\psi)+Vg^2|\psi|^2\Big].
\end{equation}
Observe that
\begin{multline}\label{eq.exp3}
	A\nabla(g\psi)\nabla(g\psi)=\psi^2 A\nabla g\nabla g+g^2A\nabla\psi\nabla\psi+\psi gA\nabla g\nabla\psi+\psi gA\nabla\psi\nabla g\\ \leq C\big[g^2A\nabla\psi\nabla\psi+ \psi^2A\nabla g\nabla g\big].
\end{multline}
Now, since $g^2A\nabla\psi\nabla\psi=A\nabla\psi\nabla(g^2\psi)-2g\psi A\nabla\psi\nabla g$, then
\begin{equation}\label{eq.exp4}
	g^2A\nabla\psi\nabla\psi\leq  2A\nabla\psi\nabla(g^2\psi)+C|\psi|^2A\nabla g\nabla g.
\end{equation}
On the other hand,   $g^2\psi\in\f D(\Omega)$, and we may write
\begin{equation}\label{eq.exp2}
\int_\Omega\Big[A\nabla\psi\nabla(g^2\psi)+Vg^2|\psi|^2\Big]=\int_\Omega g^2f\psi.
\end{equation}
Putting (\ref{eq.exp1}), (\ref{eq.exp3}), (\ref{eq.exp4}), and (\ref{eq.exp2}) together, it follows that
\begin{equation}\label{eq.exp5}
\int_\Omega\frac1ug^2|\psi|^2\leq C\Big[\int_\Omega\psi^2  A\nabla g\nabla g+\int_\Omega g^2|f\psi|\Big],
\end{equation}
where $C$ depends only on $\lambda$. From (\ref{eq.exp5}) and Lemma \ref{lm.lip}, we see that
\[
\int_\Omega\frac1ue^{2\ep\phi_M}\Big[1-\ep C\Big]\psi^2\leq C\int_{\supp f}e^{2\ep\phi_M}uf^2+\frac14\int_\Omega\frac1ue^{2\ep\phi_M}\psi^2,
\]
whence, if $\ep=1/(2C)$, we obtain for each $M\in\bb N$ the estimate
\[
\int_\Omega\frac1ue^{2\ep\phi_M}\psi^2\leq C\int_{\Omega}uf^2
\]
Then the desired estimate (\ref{eq.lmest}) for bounded $V$ follows by Fatou's lemma.

We now consider the general case $V\in L^1_{\loc}(\Omega)$ for (\ref{eq.lmest}). Let $V_N:=\min\{V,N\}$, and $u_N$ the landscape function for $L_{V_N}$. Using Lemma \ref{lm.approx}, we have that $u_N\ra u$ strongly in $L^2(\supp f)$, $u_N\searrow u$ pointwise as $N\ra\infty$ and, under assumption (\ref{eq.nosing}), it can be shown that for each $x\in\Omega$, $\rho(x,\supp f,\frac1{u_N})\ra\rho(x,\supp f,\frac1{\hat u})$ as $N\ra\infty$. Then  (\ref{eq.lmest}) is obtained by Fatou's lemma.

The proof of  (\ref{eq.lmests}) when $A$ is symmetric and continuous is similar to the argument above, except   no dependence on $\lambda$ need be introduced. Following the argument as above, the required convergence as $M\ra\infty$, and later as $N\ra\infty$, of the first term on the left-hand side (\ref{eq.lmests}) follows by applications of Fatou's lemma, Cauchy inequality with $\ep>0$, the Dominated Convergence Theorem, and (\ref{eq.unprinciplesym}). We omit the details.\hfill{$\square$}

\begin{remark} In (\ref{eq.lmests}), we may take $\alpha\in(0,1)$   either if $\Omega$ is bounded, or if we drop the first integral on the left-hand side of (\ref{eq.lmests}).
\end{remark}

We now prove the exponential decay estimates of the integral kernels.

\noindent\emph{Proof of Theorem \ref{thm.greenexp}.} We first prove (\ref{eq.gineqavg}). Fix $x_0,y_0\in\Omega$ and let $\delta,\ep$ be positive numbers, to be chosen later, with the constraint that $B(x_0,\delta)\cup B(y_0,\ep)\subset\Omega$. Let $f:=\1_{B(y_0,\ep)}$. Then $f\in L^{\infty}_c(\Omega)$, and therefore the identity
\[
(L^{-1}f)(x)=\int_{\Omega}G(x,z)f(z)\,dz=\int_{B(y_0,\ep)}G(x,y)\,dz
\]
holds for a.e. $x\in\Omega$. Then we may use (\ref{eq.lmest}) to see that
\begin{equation}\label{eq.greendecay1}\nonumber
	\int_{\Omega}\frac{e^{2\alpha\rho(x,B(y_0,\ep))}}{u(x)}\Big|\int_{B(y_0,\ep)}G(x,z)\,dz\Big|^2\,dx\leq\frac1{(1-\alpha^2)^2}\int_{B(y_0,\ep)}u(y)\,dy.
\end{equation}
Write $\rho(x)=\rho(x,B(y_0,\ep))$. We have that
\begin{multline}\label{eq.clevercs}\nonumber
\Big|\int_{B(x_0,\delta)}\int_{B(y_0,\ep)}G(x,z)\,dz\,dx\Big|^2\\=\Big|\int_{B(x_0,\delta)}\Big(\int_{B(y_0,\ep)}G(x,z)\,dz\Big)\frac{e^{\alpha\rho(x)}}{\sqrt{u(x)}}\sqrt{u(x)}e^{-\alpha\rho(x)}\,dx\Big|^2\\[3mm]\leq\Big\{\int_{B(x_0,\delta)}\frac{e^{2\alpha\rho(x)}}{u(x)}\Big|\int_{B(y_0,\ep)}G(x,z)\,dz\Big|^2\,dx\Big\}\Big\{\int_{B(x_0,\delta)}u(x)e^{-2\alpha\rho(x)}\,dx\Big\}\\[3mm]\leq\Big\{\frac1{(1-\alpha^2)^2}\int_{B(y_0,\ep)}u(y)\,dy\Big\}\Big\{\int_{B(x_0,\delta)}u(x)e^{-2\alpha\rho(x)}\,dx\Big\},
\end{multline}
where we used the Cauchy-Schwartz inequality and (\ref{eq.lmest}). Hence,
\begin{equation}\label{eq.gineq}\nonumber
\int_{B(x_0,\delta)}\int_{B(y_0,\ep)}G(x,y)\,dy\,dx\leq\frac1{(1-\alpha^2)}\Big(\int_{B(y_0,\ep)}u\Big)^{\frac12}\Big(\int_{B(x_0,\delta)}ue^{-2\alpha\rho}\Big)^{\frac12}, 
\end{equation}
or in other words,
\begin{equation} \nonumber
	\dashint_{B(x_0,\delta)}\dashint_{B(y_0,\ep)}G(x,y)\,dy\,dx\leq\frac1{|B(0,1)|}\frac1{(1-\alpha^2)}\delta^{-\frac n2}\ep^{-\frac n2}\Big(\dashint_{B(y_0,\ep)}u\Big)^{\frac12}\Big(\dashint_{B(x_0,\delta)}ue^{-2\alpha\rho}\Big)^{\frac12},
\end{equation}
which finishes the proof of (\ref{eq.gineqavg}).

We now prove (\ref{eq.greendecay}). Thus we take $n\geq3$, $\delta=\ep=\Vert V\Vert_{L^{\infty}(\Omega)}^{-1/2}$ and assume that $B(x_0,4\delta)\cap B(y_0,4\ep)=\varnothing$, $B(x_0,4\delta)\cup B(y_0,4\ep)\subset\Omega$. In this case, the Green's function $G$ is known (see \cite{dhm}, for instance) to be a continuous   weak solution of the equation $LG(x,\cdot)=0$ away from $x$, and moreover, it verifies that $G(x,y)=G^T(y,x)$ for each $x,y\in\Omega$, $x\neq y$. Let $\tilde x$, $\tilde y$ be the points associated to the balls $B(x_0,\delta), B(y_0,\ep)$ from Proposition \ref{prop.harnacku}. By the Harnack inequality \cite[Theorems 8.17 and 8.18]{gtpde} applied for $G(\cdot,y_0)$ and then for $G^T(\tilde x,\cdot)$, it is easy to see that $G(x_0,y_0)\approx G(\tilde x,\tilde y)$.  Similarly, we obtain that
\[
G(\tilde x,\tilde y)\lesssim\dashint_{B(\tilde x,c_2\delta)}\dashint_{B(\tilde y,c_2\ep)}G(x,y)\,dy\,dx,
\]
where $c_2$ is the constant from Proposition \ref{prop.harnacku}.  By Proposition \ref{prop.harnacku} again, we have that
\[
\dashint_{B(\tilde y,c_2\ep)}u\lesssim  u(\tilde y),\qquad \dashint_{B(\tilde x,c_2\delta)}u\lesssim u(\tilde x).
\]
By the Triangle inequality and Proposition \ref{prop.harnacku}, for all $x\in B(\tilde x,c_2\delta)$ we have that
\[
\rho\Big(x,B(y_0,\ep),\frac1u\Big)\geq\rho\Big(\tilde x,\tilde y,\frac1u\Big)-C.
\]
Therefore, (\ref{eq.gineqavg}) implies that
\begin{equation}\label{eq.withconst}
G(\tilde x,\tilde y)\lesssim \sqrt{u(\tilde y)}\sqrt{u(\tilde x)}\Vert V\Vert_{L^{\infty}(\Omega)}^{\frac n2}e^{-\alpha\rho(\tilde x,\tilde y,\frac1u)},
\end{equation}
which gives the desired result.\hfill{$\square$}

\begin{remark}\label{rm.neumann} Under the setting of Theorem \ref{thm.greenexp}, if $V\in L^{\infty}(\Omega)$ and $\Omega$ is bounded, then we may define a Neumann landscape function $u_{\f N}$ on $\Omega$, via the operator with 0 Neumann boundary conditions, similarly as was done in Section \ref{sec.optheory}. Using this landscape function, it is possible to obtain the following analogue of (\ref{eq.greendecay}): 
\begin{equation}\label{eq.greendecay2}
	G_{\Omega}(x,y)\lesssim\Vert V\Vert_{L^{\infty}(\Omega)}^{\frac n2}\sqrt{u_{\f N}(x)}\sqrt{u_{\f N}(y)}e^{-\alpha\rho(x,y,\frac1{u_{\f N}})},
\end{equation} 
where $G_{\Omega}$ is the integral kernel for the (Dirichlet) solution operator $L^{-1}$ on $\Omega$. The proof of (\ref{eq.greendecay2}) is similar to (and simpler than)   that of (\ref{eq.greendecay}).	
\end{remark}

\begin{remark}\label{rm.expdecay} Take the assumptions of Theorem \ref{thm.greenexp}, $n\geq3$, and suppose it were known that $U:=\Vert u\Vert_{L^{\infty}(\Omega)}<\infty$.  Then we may use (\ref{eq.gineqavg}) with $\delta=\ep=|x_0-y_0|/3$ together with the fact that $\frac1u\geq\frac1U$ and the Moser estimate to obtain that
\[
G(x_0,y_0)\lesssim\frac1{|x_0-y_0|^n}Ue^{-\frac{\alpha}3\frac{|x_0-y_0|}{\sqrt U}},
\]
whenever $B(x_0,4\delta)\cup B(y_0,4\ep)\subset\Omega$. Now, if $x_0,y_0$ are such that $|x_0-y_0|^2>U$, then 
\[
G(x_0,y_0)\lesssim\frac1{|x_0-y_0|^{n-2}}e^{-\frac{\alpha}3\frac{|x_0-y_0|}{\sqrt U}}.
\]
Otherwise, if $|x_0-y_0|^2\leq U$, then we have by the domination of the Green's function by the fundamental solution $\Gamma(x_0,y_0)$ that
\[
G(x_0,y_0)\lesssim\frac1{|x_0-y_0|^{n-2}}e^{-\frac{\alpha}3\frac{|x_0-y_0|}{\sqrt U}}.
\]
Hence, in any case, we obtain that
\begin{equation}\label{eq.expdecay}
	G(x_0,y_0)\lesssim\frac1{|x_0-y_0|^{n-2}}e^{-\frac{\alpha}3\frac{|x_0-y_0|}{\sqrt U}},
\end{equation}
for each $x_0,y_0\in\Omega$ such that $B(x_0,\tfrac43|x_0-y_0|)\cup B(y_0,\tfrac43|x_0-y_0|)\subset\Omega$. The implicit constant in (\ref{eq.expdecay}) does not depend on $u$, $V$, $U$ nor on $\Omega$.
\end{remark}

In order to study the exponential decay of eigenfunctions, we must first consider a non-homogeneous analogue of the   operator $L$ defined in Section \ref{sec.prelim}. Recall that $\f D(\Omega)$ is the completion of $C_c^{\infty}(\Omega)$ in the norm $\Vert\psi\Vert_{\f l}$ given by (\ref{eq.norm}). 

\begin{definition}[The non-homogeneous operator $\n L$]\label{def.nonhomo} We let $\n D(\Omega)$ be the completion of $C_c^{\infty}(\Omega)$ in the norm
\[
\sqrt{\Vert\psi\Vert_{\f l}^2+\Vert\psi\Vert_{L^2(\Omega)}^2},\qquad\psi\in C_c^{\infty}(\Omega).
\]
Following Ouhabaz \cite{ou} (see also Section 2.1 of \cite{mp} for more details), we may associate to the form $\f l$ an unbounded operator $\n L:\n D(\n L)\ra L^2(\Omega)$ where  
\begin{multline}\nonumber
\n D(\n L)=\Big\{\psi\in\n D(\Omega) \text{ such that there exists }  f\in L^2(\Omega)\\ \text{ with } \f l(\psi,\varphi)=\langle f,\varphi\rangle_{L^2(\Omega)}   \text{ for every } \varphi\in \n D(\Omega)\Big\},\qquad \n L\psi:=f.
\end{multline}
Then Proposition 1.22 in \cite{ou} applies and we conclude that $\n L$ is densely defined, for every $\ep>0$ the operator $\n L+\ep$ is invertible from $\n D(\n L)$ into $L^2(\Omega)$, and its inverse $(\n L+\ep)^{-1}$ is a bounded operator on $L^2(\Omega)$. In addition,
\[
\Vert\ep(\n L+\ep)^{-1}f\Vert_{L^2(\Omega)}\leq\Vert f\Vert_{L^2(\Omega)},\qquad\text{for each }\ep>0, f\in L^2(\Omega).
\]
\end{definition}

For each $t\in(0,\infty)$, and $f\in L^2(\Omega)$,   denote
\[
V_t=V+\frac1{t^2},\quad \n L_t=\n L+\frac1{t^2},\quad \n R_tf:=(1+t^2\n L)^{-1}f=\frac1{t^2}\n L_{t}^{-1}f
\]
We call  $\n R_t$ a \emph{resolvent of $\n L$}. Observe that   $L_t(\n R_tf)=\frac1{t^2}f$ in the weak sense. Let $u_t$   be the landscape function of Theorem \ref{thm.exist} for the operator $L_t=-\dv A\nabla+(V+\frac1{t^2})$. The following result is not, strictly speaking, a corollary of Theorem \ref{thm.lmr}, but its proof is essentially the same.

\begin{corollary}[Exponential decay estimate for resolvents]\label{cor.resexp} Retain the setting and assumptions of Theorem \ref{thm.lmr}.  	  For each $t\in(0,\infty)$, let $u_t$ and $\n R_tf$ be as defined above. Then for all $\alpha\in(0,1)$ and $f\in L^2(\Omega)$ with compact support, we have the estimate
	\begin{equation}\label{eq.resest2} 
		\int_{\Omega}\frac1{u_t}e^{2\alpha\rho(\cdot,\supp f,\frac1{\hat u_t})}|\n R_tf|^2~\leq C\frac1{t^4}\int_{\Omega} u_tf^2,
	\end{equation}
	where $C$ depends only on $\lambda$ and $\alpha$. Moreover, if $A$ is symmetric and continuous, we have a corresponding analogue of (\ref{eq.lmests}).
\end{corollary}

\begin{remark} Note that, for each $t>0$, $u_t$ can be replaced by $u$ in (\ref{eq.resest2}).
\end{remark}

\noindent\emph{Proof of Theorem \ref{thm.eigen}.} If $A$ is symmetric, then the operator $\n L$ from Definition \ref{def.nonhomo} is self-adjoint; moreover, simple coercivity arguments show that $\n L$ has no negative eigenvalues. The proof of the estimate (\ref{eq.eigen}) is essentially the same as that of (\ref{eq.lmests}) and (\ref{eq.lmest}) given above; one would need to use that $\n D(\n L)\subset\n D(\Omega)\subset\f D(\Omega)$ and that $\rho(x,E,w)=0$ for each $x\in E$. We omit further details.\hfill{$\square$}

\section{Comparison to the Fefferman-Phong-Shen maximal function}\label{sec.fps}

We will  heavily make use of the following properties of the maximal function $m$.

\begin{proposition}[Properties of the Fefferman-Phong-Shen maximal function \cite{shenf}] Retain the settings and assumptions of Theorem \ref{thm.compare}. Then
	\begin{enumerate}[(a)]
		\item $0<m(x)<\infty$ for every $x\in\bb R^n$,\
		\item if $y\in B(x, C/m(x))$, then $m(y)\approx m(x)$.
		\item there exists $k_0>0$ such that for any $x,y\in\bb R^n$,
		\[
		m(y)\gtrsim\frac{m(x)}{\big[1+|x-y|m(x)\big]^{\frac{k_0}{k_0+1}}}.
		\]
	\end{enumerate}	
\end{proposition}

\noindent\emph{Proof of Theorem \ref{thm.compare}.} First, since $V$ is a Shen potential, then (\ref{eq.greenint}) is known to hold \cite{shenf, mp}, so that the function $u$ from Theorem \ref{thm.exist} solves $Lu=1$ weakly in $\bb R^n$ and $u=\int_{\bb R^n}\Gamma(\cdot,y)\,dy$, where $\Gamma$ is the fundamental solution for the operator $L$.  Now we  show the lower bound for $u$. Fix $x\in\bb R^n$. By \cite[Lemma 4.8]{shenf}\footnote{See also \cite[Lemma 7.17]{mp} for the case of $-\dv A\nabla+V$.}, we have that 
\[
\Gamma(x,y)\geq\frac{c_n}{2|x-y|^{n-2}},\qquad\text{for each }y\in B\Big(x,\frac c{m(x)}\Big),
\]
where $c_n$ is the dimensional constant such that $\Gamma_{-\Delta}=\frac{c_n}{|x-y|^{n-2}}$. It follows that
\begin{equation}\nonumber
	u(x)=\int_{\bb R^n}\Gamma(x,y)\,dy\geq\int_{B(x,c/m(x))}\Gamma(x,y)\,dy\gtrsim\int_0^{c/m(x)}\frac1{r^{n-2}}r^{n-1}\,dr  \gtrsim \frac1{m^2(x)}.
\end{equation}

We turn to the proof of the upper bound for $u$, which is more involved. First, note that by \cite[Theorem 4.16]{mp}, we may procur the estimate 
\begin{equation}\label{eq.estm5}
	\int_{\bb R^n}m^2(y)e^{2\ep\rho(y,\supp f,m^2)}|L^{-1}f(y)|^2\,dy \lesssim\int_{\bb R^n}\frac1{m^2(y)}f^2(y)\,dy
\end{equation}
for each $f\in L^2(\bb R^n)$ with compact support\footnote{We may also prove (\ref{eq.estm5}) in a similar way as we proved (\ref{eq.lmest}).}. Note that the landscape function $u$ over $\bb R^n$ is not a Lax-Milgram solution; however, we may write it as a sum of Lax-Milgram solutions by decomposing the constant function $1$ into annuli; then we will exploit the fact that these decompositions exhibit exponential decay by virtue of (\ref{eq.estm5}) to finish the proof. Fix $x\in\bb R^n$, let $A_0:= B\big(x,\frac1{m(x)}\big)$, and for each $k\in\bb N$  let $A_k:=B\big(x,\frac{2^k}{m(x)}\big)\backslash B\big(x,\frac{2^{k-1}}{m(x)}\big)$. Write 
\[
f_k:=\1_{A_k}, \qquad u_k:= L^{-1}f_k,\qquad U_M:=\sum_{k=0}^Mu_k.
\]
By the maximum principle, we have that $U_N\geq U_M$ whenever $N\geq M$. Moreover, since
\[
u(z)=\int_{\bb R^n}\Gamma(z,y)\,dy=\lim_{M\ra\infty}\int_{\bb R^n}\Gamma(z,y)\sum_{k=0}^M\1_{A_k}(y)\,dy=\lim_{M\ra\infty} U_M,
\]
it follows that $U_M\nearrow u$ pointwise as $M\ra\infty$.

We now estimate each $u_k(x)$. For $u_0(x)$, by (\ref{eq.estm5}) we have that
\begin{equation}\nonumber
\int_{B(x,1/m(x))}u_0^2(y)\,dy\lesssim\int_{B(x,1/m(x))}\frac1{m^4(y)}\,dy\lesssim\frac1{m^4(x)}|B(x,1/m(x))|,
\end{equation}
where we have also used the slowly-varying properties of $m$.  Thus we have that
\[
\Big(\dashint_{B(x,1/m(x))}u_0^2(y)\,dy\Big)^{1/2}\lesssim\frac1{m^2(x)}.
\]
Then, applying the Moser estimate for   $-\dv A\nabla w+Vw=1$, we conclude that
\begin{equation}\label{eq.u0}
u_0(x)\lesssim\frac1{m^2(x)}.
\end{equation}
Now let $k\geq1$. Using (\ref{eq.estm5}), we see that
\begin{multline}\nonumber
\int_{B(x,1/m(x))}u_k^2(y)e^{2\ep\rho(y,\supp f_k, m^2)}\,dy\lesssim\frac1{m^2(x)}\int_{A_k}\frac1{m^2(y)}\,dy\\ \lesssim\frac1{m^4(x)}\int_{A_k}\big[1+|y-x|m(x)\big]^{\frac{a}{a+1}}\,dy\lesssim\frac1{m^4(x)}\int_{2^{k-1}/m(x)}^{2^k/m(x)}\big[1+rm(x)\big]^{\frac{a}{a+1}}r^{n-1}\,dr\\ \lesssim\frac1{m^4(x)}|B(x,1/m(x))|2^{2kn},
\end{multline}
where we have used the properties of $m$. Using these properties once again and the Moser estimate (note that $u_k$ solves $-\dv A\nabla w+Vw=0$ in $B(x,1/m(x))$), we obtain that
\begin{multline}\label{eq.ukest}
u_k(x)\lesssim\Big(\dashint_{B(x,1/m(x))}u_k^2(y)\,dy\Big)^{1/2}\lesssim\frac1{m^2(x)}2^{kn}e^{-\ep\rho(x,\supp f_k,m^2)}\\ \lesssim\frac1{m^2(x)}2^{kn}e^{-\ep[1+\dist(x,\supp f)m(x)]^{\frac1{a+1}}}\lesssim\frac1{m^2(x)}2^{kn}e^{-\ep2^{\frac k{a+1}}}.
\end{multline}

We are ready to collect our estimates. For each $M\in\bb N$, use (\ref{eq.u0}) and (\ref{eq.ukest}) to see that
\[
U_M(x)=\sum_{k=0}^Mu_k(x)\lesssim\frac1{m^2(x)}\sum_{k=0}^M2^{kn}e^{-\ep2^{\frac k{a+1}}}\lesssim\frac1{m^2(x)}.
\]
Since $U_M(x)\nearrow u(x)$ and $U_M$ is  bounded by $1/m^2(x)$, the desired result follows.\hfill{$\square$}

\section{Applications to the magnetic Schr\"odinger operator}\label{sec.ms}

\subsection{The uncertainty principle for $L_{\a,V}$ and a priori exponential decay results}\label{sec.expms} First, let us prove the uncertainty principle for the magnetic Schr\"odinger operator via the landscape function for the operator $\tilde L$ from (\ref{eq.electric}).

\noindent\emph{Proof of Theorem \ref{thm.magup}.}  By assumption, we have that (perhaps after an orthonormal change of variables) $\Sigma_{\m S}\B+V\geq0$ on $\Omega$. By the uncertainty principle (\ref{eq.unprinciplesym}), we may write, for any real-valued Lipschitz function $g$ on $\Omega$ with compact support,
\begin{equation}\label{eq.fm}
	\int_\Omega u^2|\nabla(g/u)|^2+\int_{\Omega}\frac1ug^2\leq\int_{\Omega}|\nabla g|^2+(\Sigma_{\m S}\B+V)g^2.
\end{equation}
It follows that, for any complex-valued $f\in C_c^{\infty}(\bb R^n)$, by taking $g=|f|$,
\begin{equation}\label{eq.fm2}
	\int_{\Omega}\frac1u|f|^2\leq\int_{\Omega}\big|\nabla|f|\big|^2+(\Sigma_{\m S}\B+V)|f|^2\leq\int_{\Omega}|D_{\bf a}f|^2+(\Sigma_{\m S}\B+V)|f|^2,
\end{equation}
where in the last inequality we used the diamagnetic inequality $|\nabla|f||\leq|D_{\bf a}f|$.

Now, define the commutator
\begin{equation}\label{eq.com}
	[D_j,D_k]=D_jD_k-D_kD_j,
\end{equation}
and note that
\begin{equation}\label{eq.field}
	b_{jk}(x)f(x)=-i([D_j,D_k]f)(x),\qquad\text{for each }x\in\Omega.
\end{equation}
Then, we may estimate the middle term of (\ref{eq.fm2}) as follows:
\begin{align*} 
\Big|\int_{\Omega}\Sigma_{\m S}\B(x)|f(x)|^2\,dx\Big|&=\Big| \sum_{(j',k')\in\m S}\int_{\Omega}b_{j'k'}(x)f(x)\overline{f(x)}\,dx\Big|\\&= \Big|\sum_{(j',k')\in\m S}\int_{\Omega}\big\{(D_{j'}D_{k'}-D_{k'}D_{j'})f\big\}\overline{f(x)}\,dx\Big|\\ &= \Big|\sum_{(j',k')\in\m S}\int_{\Omega}\Big[-(D_{j'}D_{k'}f)(x)\overline{f(x)}+(D_{k'}D_{j'}f)(x)\overline{f(x)}\Big]\,dx\Big|\\ &=  \Big|\sum_{(j',k')\in\m S}\int_{\Omega}\big[-D_{k'}f\overline{D_{j'}f})+D_{j'}f\overline{D_{k'}f}\big]\Big|\\ &\leq  \sum_{(j',k')\in\m S}\int_{\Omega}|D_{k'}f|^2+|D_{j'}f|^2\\ &\qquad\qquad=(n-1)\int_{\Omega}\sum_{k=1}^n|D_kf|^2=(n-1)\int_{\Omega}|D_{\bf a}f|^2.
\end{align*}
The desired result follows.\hfill{$\square$}

\subsubsection*{Operator-theoretic preliminaries} Now, let us give the necessary setup for Theorem \ref{thm.expdecayms}; see Section 2.1 of \cite{mp} for a more general setting, or more details. For the rest of this section, all functions are assumed complex-valued, unless explicitly mentioned otherwise. On $C_c^{\infty}(\Omega)$ define the form
\begin{equation}\label{eq.aform}
\f a(\psi,\varphi)=\int_\Omega\Big[D_\a\psi\overline{D_\a\varphi}+V\psi\overline{\varphi}\Big],\qquad\psi,\varphi\in C_c^{\infty}(\Omega).
\end{equation}
Under the assumptions of Theorem \ref{thm.expdecayms}, we may prove, similar to Theorem \ref{thm.magup}, that
\begin{equation}\label{eq.vbound}
\int_{\Omega}V_-|\psi|^2\leq\frac{n-1}n\int_\Omega|D_\a\psi|^2,\quad\text{for each }\psi\in C_c^{\infty}(\Omega).
\end{equation}
It follows that
\begin{equation}\label{eq.vbound2}
\int_\Omega V_-|\psi|^2\leq(n-1)\f a(\psi,\psi),\quad\text{for each }\psi\in C_c^{\infty}(\Omega).
\end{equation}
From (\ref{eq.aform}), (\ref{eq.vbound}), (\ref{eq.vbound2}), the diamagnetic inequality, and the fact that $\a$ and $V$ are not both identically $0$, we have that $\f a(\psi,\psi)\geq0$, and $\f a(\psi,\psi)=0$ implies  $\psi\equiv0$. Now, using (\ref{eq.vbound}), (\ref{eq.vbound2}), and the Cauchy-Schwartz inequality, it is easily seen that
\[
|\f a(\psi,\varphi)|\leq9n\Vert\psi\Vert_{\f a}\Vert\varphi\Vert_{\f a},
\]
where
\[
\Vert\psi\Vert_{\f a}:=\sqrt{\f a(\psi,\psi)}.
\]
Moreover, from (\ref{eq.up}) and (\ref{eq.vbound2}), we have that
\begin{equation}\label{eq.magup1}
\int_\Omega u_{\m S}^2\Big|\nabla\Big(\frac{|\psi|}{u_{\m S}}\Big)\Big|^2+	\int_\Omega\frac1{u_{\m S}}|\psi|^2\leq2n(n-1)\int_{\Omega}\Big[|D_{\bf a}\psi|^2+V|\psi|^2\Big],\qquad\psi\in C_c^{\infty}(\Omega).
\end{equation}
We may now set up Lax-Milgram solutions. Since the form $\f a$ is symmetric and sesquilinear, we know that $\f a$ is an inner product on $C_c^{\infty}(\Omega)$, and therefore $\Vert\cdot\Vert_{\f a}$ is a norm on $C_c^{\infty}(\Omega)$.  Now we let $\f D^{\f a}(\Omega)=\f D^{\f a}_{\a,V}(\Omega)$ be the completion of $C_c^{\infty}(\Omega)$ in the norm $\Vert\cdot\Vert_{\f a}$. Then $(\f D^{\f a}(\Omega),\f a)$ is a Hilbert space on which $\f a$ is a bounded and coercive form. Then (\ref{eq.vbound}), (\ref{eq.vbound2}), and (\ref{eq.magup1}) hold  also for $\psi\in\f D^{\f a}(\Omega)$.  Let $L_{\a, V}:\f D^{\f a}(\Omega)\ra(\f D^{\f a}(\Omega))^*$ be the operator defined as follows: for $\psi\in\f D^{\f a}(\Omega)$, $L_{\a, V}\psi$ is the functional given by the rule $\langle L_{\a, V}\psi,\varphi\rangle=\f a(\psi,\varphi)$ for any $\varphi\in\f D^{\f a}(\Omega)$. Then, by the Lax-Milgram theorem, $L_{\a, V}$ is invertible. Via a compactness argument similar to the proof of (\ref{eq.coercive}), it can be shown that $L^q_c(\Omega)\subset(\f D^{\f a}(\Omega))^*$ for $q>\max\{1,\frac{2n}{n+2}\}$. For $n\geq3$, the existence of Green's function (as an integral kernel) for the operator $L_{\a, V}$ under our assumptions was shown in \cite[Theorem 5.17]{mp}\footnote{The proof there was for the fundamental solution; but it can be adapted to Lipschitz domains.}.

To set up eigenfunctions of the magnetic Schr\"odinger operator, we consider a non-homogeneous version of $L_{\a, V}$, similar to Definition \ref{def.nonhomo}. Let $\n D^{\f a}(\Omega)$ be the completion of $C_c^{\infty}(\Omega)$ in the norm $\sqrt{\f a(\psi,\psi)+\Vert\psi\Vert_{L^2(\Omega)}^2}$. Then, in the terminology of \cite{ou}, $\f a$ is a densely defined, accretive, continuous, and closed sesquilinear form on $L^2(\Omega)$. We may thus associate to $\f a$ an unbounded operator $\n L^{\f a}:\n D(\n L^{\f a})\ra L^2(\Omega)$ with
\begin{multline}\nonumber
	\n D(\n L^{\f a})=\Big\{\psi\in\n D^{\f a}(\Omega) \text{ such that there exists }  f\in L^2(\Omega)\\ \text{ with } \f a(\psi,\varphi)=\langle f,\varphi\rangle_{L^2(\Omega)}   \text{ for every } \varphi\in \n D^{\f a}(\Omega)\Big\},\qquad \n L^{\f a}\psi:=f.
\end{multline}
Note that $\n D(\n L^{\f a})\subset\n D^{\f a}(\Omega)\subset L^2(\Omega)$. Then from Proposition 1.22   \cite{ou}  we conclude that $\n L^{\f a}$ is densely defined, for every $\ep>0$ the operator $\n L^{\f a}+\ep$ is invertible from $\n D(\n L^{\f a})$ into $L^2(\Omega)$, and its inverse $(\n L^{\f a}+\ep)^{-1}$ is a bounded operator on $L^2(\Omega)$. In addition,
\[
\Vert\ep(\n L^{\f a}+\ep)^{-1}f\Vert_{L^2(\Omega)}\leq\Vert f\Vert_{L^2(\Omega)},\qquad\text{for each }\ep>0, f\in L^2(\Omega).
\]
In particular, $\n L^{\f a}$ is self-adjoint, and has no negative eigenvalues.

\noindent\emph{Proof of Theorem \ref{thm.expdecayms}.} We first prove \ref{item.lmms}; the technique is similar to that of Theorem \ref{thm.lmr}, but   the possible   negative part of $V$ is a bit cumbersome. We show how to obtain an analogue estimate to (\ref{eq.exp5}) when $\Sigma_S\B+V$ is bounded, from which the desired result in the full generality will follow in the same way as it did in the proof of (\ref{eq.lmest}). Fix $f\in L^2_c(\Omega)$, and $\psi:=L_{\a,V}^{-1}f$.  If $g=e^{\phi}$ where $\phi$ is Lipschitz and bounded, then $g\psi,g^2\psi\in\f D^{\f a}(\Omega)$. By the Cauchy inequality with $\ep>0$,
\begin{equation}\label{eq.id}
|D_{\a}(g\psi)|^2\leq(1+\ep)g^2|D_{\a}\psi|^2+C_\ep|\nabla g|^2|\psi|^2.
\end{equation}
On the other hand,  $D_{\a}\psi\overline{D_{\a}(g^2\psi)}=g^2|D_{\a}\psi|^2+2g\overline{\psi}\nabla g D_{\a}\psi$, so in turn we have
\begin{equation}\label{eq.id2}
g^2|D_{\a}\psi|^2\leq\frac1{1-\ep}\Re e\big\{D_{\a}\psi\overline{D_{\a}(g^2\psi)}\big\}+C_\ep|\nabla g|^2|\psi|^2.
\end{equation}
Moreover, since $L_{\a,V}\psi=f$, then
\begin{equation}\label{eq.id3}
\int_\Omega D_{\a}\psi\overline{D_{\a}(g^2\psi)}=\int_\Omega fg^2\overline{\psi}-\int_\Omega Vg^2|\psi|^2.
\end{equation}
Putting together (\ref{eq.id}), (\ref{eq.id2}), and (\ref{eq.id3}), we see that
\begin{equation}\label{eq.cauchy}
\int_\Omega|D_{\a}(g\psi)|^2\leq\frac{1+\ep}{1-\ep}\int_\Omega V_-g^2|\psi|^2-\frac{1+\ep}{1-\ep}\int_\Omega V_+g^2|\psi|^2  + C_\ep\Big[\int_\Omega g^2|f||\psi|+\int_\Omega|\nabla g|^2|\psi|^2\Big].
\end{equation}
Using (\ref{eq.vbound}) now, we may choose $\ep$ small (depending only on $n$) so that the first term on the right-hand side of (\ref{eq.cauchy}) can be absorbed to the left-hand side. Consequently, from this estimate and the uncertainty principle (\ref{eq.up}), we obtain that
\begin{equation}\label{eq.magup2}
	\int_\Omega u_{\m S}^2\Big|\nabla\Big(\frac{g|\psi|}{u_{\m S}}\Big)\Big|^2+	\int_\Omega\frac1{u_{\m S}}g^2|\psi|^2\leq C_n \Big[\int_\Omega g^2|f||\psi|+\int_\Omega|\nabla g|^2|\psi|^2\Big].
\end{equation}
From (\ref{eq.magup2}), it is not hard to conclude the proof of \ref{item.lmms}.

The proof of \ref{item.eigenms} follows similarly, using the estimate (\ref{eq.magup2}); we skip the details.

We turn to \ref{item.greenms}. Since $G_{\a,V}$ is in general complex-valued, the method of proof of (\ref{eq.greendecay}) will not work here. Instead, we follow a method of proof similar to Shen's upper bound exponential decay estimate \cite{shenf}; see also \cite[Section 6]{mp}. To start, fix a ball $B\subset\Omega$, let $\varphi\in C_c^{\infty}(\Omega\backslash 2B;\bb R)$ with $\varphi\equiv0$ on $2B$, and note that if $\psi\in L^2_{\loc}(\Omega\backslash B)$  and $D_\a\psi\in L^2_{\loc}(\Omega\backslash B)$ with $L_{\a,V}\psi=0$ in the weak sense on $\Omega\backslash B$, then using (\ref{eq.magup2}), we may prove that
\begin{equation}\label{eq.profo}
\int_\Omega\frac1{u_{\m S}}\varphi^2|\psi|^2g^2\leq C_n\int_\Omega g^2|\nabla\varphi|^2|\psi|^2,
\end{equation}
where $g=e^{\ep\phi}$, and $\phi$ is Lipschitz, real, and bounded on $\Omega$. Fix $x,y,B_x,B_y$ as in the statement of \ref{item.greenms}, and let $\tilde x,\tilde y$ be the points associated to the balls $B_x$, $B_y$ from Proposition \ref{prop.harnacku}. Note that $L_{\a, V}G_{\a, V}(\cdot,\tilde y)=0$ in the weak sense in $\Omega\backslash B(\tilde y,\frac{c_2}4\B_\infty^{-1/2})=\Omega\backslash\tilde B$, where $c_2$ is the constant from Proposition \ref{prop.harnacku}. Fix $M\in\bb N$ large and choose $\varphi\in C_c^{\infty}(B(\tilde y,2M))$ as the cut-off function $0\leq\varphi\leq1$, $\varphi\equiv0$ on $2\tilde B$, $|\nabla\varphi|\lesssim2\B_\infty^{1/2}$ on $4\tilde B\backslash2\tilde B$, $\varphi\equiv1$ on $B(\tilde y,M)\backslash4\tilde B$, and $|\nabla\varphi|\lesssim2/M$ on $B(\tilde y,2M)\backslash B(\tilde y,M)$. Use (\ref{eq.profo}) with $\psi=G(\cdot,\tilde y)$ and $\varphi$ as described above; then, since $G_{\a, V}$ verifies (\ref{eq.greenbound2}), we eventually obtain
\begin{equation}\label{eq.profo2}
\int_{\hat B}\frac1{u_{\m S}}|G(\cdot,\tilde y)|^2e^{2\ep\phi}\leq C_n\B_\infty\int_{4\tilde B\backslash2\tilde B}e^{2\ep\phi}|G(\cdot,\tilde y)|^2,
\end{equation}
after passing $M\ra\infty$ (see \cite[Section 6]{mp} for more details on a similar argument), and where $\hat B:=B(\tilde x,c_2\B_\infty^{-1/2})$. Moreover, by an approximation argument as in the proof of Theorem \ref{thm.lmr}, we may now without loss of generality take $\phi=\rho(\cdot,\tilde y,\frac1{u_{\m S}})$ in (\ref{eq.profo2}). From Proposition \ref{prop.harnacku}, we have that $u_{\m S}(z)\approx u_{\m S}(\tilde x)$ for $z\in\hat B$, and $u_{\m S}(w)\approx u_{\m S}(\tilde y)$ for $w\in4\tilde B$; we also have that $\max\{\rho(z,\tilde x,\frac1{u_{\m S}}), \rho(w,\tilde y,\frac1{u_{\m S}})\}\leq C$ for $z\in\hat B, w\in4\tilde B$, respectively.  Putting all these observations together, as well as (\ref{eq.greenbound2}) and (\ref{eq.lbound}), into (\ref{eq.profo2}) yields that
\[
\Big(\dashint_{\hat B}|G(\cdot,\tilde y)|^2\Big)^{1/2}\leq C_n\B_\infty^{n/2}\sqrt{u_{\m S}(\tilde x)}\sqrt{u_{\m S}(\tilde y)}e^{-\ep\rho(\tilde x,\tilde y,\frac1{u_{\m S}})}, 
\]
which gives (\ref{eq.greendecayms}) using   scale-invariant local boundedness \cite[Definition 6.1]{mp}.\hfill{$\square$}

\begin{remark}\label{rm.neumannms} Under the assumptions of \ref{item.greenms}, if $\Omega$ is bounded then as in Remark \ref{rm.neumann} we may consider instead a Neumann landscape function $u_{\f N,\m S}$. Using this Neumann landscape function, it is possible to prove (\ref{eq.greendecayms}) directly with $x$ and $y$, without passing to $\tilde x,\tilde y$. This is because $u_{\f N,\m S}\geq\B_\infty^{-1}$ on $\Omega$.
\end{remark}

\begin{remark} It is clear that one may also get an exponential decay estimate for resolvents of $\n L^{\f a}$ similar to that of Corollary \ref{cor.resexp}.
\end{remark}

\subsection{The solutions to Shen's problems, and other corollaries}\label{sec.cor} Here, we describe some results which are almost ``free'' from the methods in the literature when combined with the  uncertainty principle Corollary \ref{cor.upm}  under the favored directionality assumptions.  This list of corollaries is not exhaustive; several other results of independent interest from the literature are also recovered under our setting of irregular magnetic fields, but we have chosen to highlight only a select few, for the sake of brevity.   Let us emphasize that the heavy lifting of the results in this section was done by the respective cited authors. In each of the following corollaries, there is an underlying landscape function $u_{\m S}$, but we choose to work exclusively with the Fefferman-Phong-Shen maximal function (\ref{eq.fpsm}), to emphasize the connection to the corresponding results in the literature; of course, due to Theorem \ref{thm.compare}, there is essentially no difference in working with $m^2$ or with $1/u$, for Shen potentials. For literature on the results related to the corollaries presented here, we kindly direct the reader to the papers cited in each corollary. 

We begin with the exponential decay of the fundamental solution for Shen potentials.

\begin{corollary}[Exponential decay of the fundamental solution; {\cite[Corollary 6.16]{mp}}]\label{cor.fund} Let $n\geq3$, and $\Omega=\bb R^n$. Assume that $\a, \B, V$ verify  (\ref{eq.bl1}) and (\ref{eq.signass}), and that Assumption \ref{ass.g} holds\footnote{As remarked in Section \ref{sec.main2a}, Assumption \ref{ass.g} holds,  in particular, if $V\geq0$ and (\ref{eq.bl1}) is true.}. Moreover, suppose that $\m S$ is an admissible selection such that $\Sigma_{\m S}\B+V$ is a non-degenerate Shen potential. Then there exist $\ep>0$ and $C\geq1$, depending only on $n$, $C_0$, $C_1$, $\delta$, and the constants from Assumption \ref{ass.g}, such that for a.e. $x,y\in\Omega$,
\begin{equation}\label{eq.fund} 
|\Gamma_{\a,V}(x,y)|\leq\frac C{|x-y|^{n-2}}e^{-\ep\rho(x,y,m^2(\cdot,\Sigma_{\m S}\B+V))}.
\end{equation}
In particular, if $V\geq0$ and (\ref{eq.ass}) holds, then there exists $\m S$ so that $|\B|=\Sigma_{\m S}\B$.
\end{corollary}

The condition (\ref{eq.gradb}) was used in \cite{mp} to prove (\ref{eq.fund})    only to verify that the uncertainty principle (\ref{eq.upm}) is true; we have seen that under the favored directionality assumption (\ref{eq.signass}), we can recover (\ref{eq.upm}) without either assumption  (\ref{eq.gradb}) or (\ref{eq.V}). The rest of the proof in \cite[Section 6]{mp} remains unchanged, giving us (\ref{eq.fund}) in our situation.

\begin{remark} We can also obtain exponential decay estimates for heat kernels \cite{kurata} under assumptions similar to those of Corollary \ref{cor.fund}, with no assumption on $\nabla\B$.
\end{remark}

In Theorem \ref{thm.expdecayms} \ref{item.eigenms}, we have  shown an a priori $L^2$ exponential decay estimate for eigenfunctions of $L_{\a,V}$. The following result  tells us that under further assumptions, exponential decay in the Agmon distance with $1/u_{\m S}$   implies exponential decay in  the Euclidean distance. The proof   is similar to those of \cite[Theorem 0.20]{shene} and  (\ref{eq.eigenms}).

\begin{corollary}[Exponential decay of eigenfunctions for $L_{\a,V}$; {\cite[Theorem 0.20]{shene}}]\label{cor.expeigen} Let $n\geq3$, and $\Omega=\bb R^n$. Assume that $\a, \B, V$ verify  (\ref{eq.bl1}) and (\ref{eq.signass}), and that $V\geq0$. Moreover, suppose that  $|\B|+V$ is a non-degenerate Shen potential, and that there exist $\mu>0$ and $\psi\in\n D(\n L^{\f a})$ with $\n L^{\f a}\psi=\mu\psi$. Let $E_{\sigma}:=\big\{x\in\bb R^n:m^2(x,|\B|+V)\leq\sigma\big\}$. Then there exist $\ep>0$ and $C\geq1$, depending only on $n$, $C_0$, $C_1$, $\delta$,  such that 
\begin{equation}\label{eq.eigenms2}\nonumber
|\psi(x)|\leq C\mu^{n/2}e^{-\ep\rho(x,E_{C\mu},m^2(\cdot,|\B|+V))}\Vert\psi\Vert_{L^2(\bb R^n)},\quad\text{for each }x\in\Omega,\qquad\text{if }\mu\geq C.
\end{equation}
\end{corollary}

Now, we turn to   bounds on the eigenvalue counting function, when $V\geq0$. For $\mu\in\bb R$, let $N(\mu,L_{\a,V})$ be the dimension of the spectral projection for $\n L^{\f a}_{\a,V}$ (see Section \ref{sec.expms}) corresponding to the interval $(-\infty,\mu)$. When $N(\mu,L_{\a,V})$ is finite for a given $\mu$, $N(\mu, L_{\a,V})$ is in particular the number of eigenvalues (counting multiplicity) of $\n L^{\f a}_{\a,V}$ smaller than $\mu$. In the case that $V\geq0$, $\n L^{\f a}_{\a, V}$ has no negative eigenvalues. Given $\mu>0$, let $\{Q_j\}_j$ be a grid of mutually disjoint cubes on $\bb R^n$ of sidelength $\frac1{\sqrt{\mu}}$, and let $\tilde N(\mu)$ be the number of cubes $Q_j$ such that 
\begin{equation}\label{eq.tilden}
	\Big(\dashint_{Q_j}|\B|^{n/2}\Big)^{2/n}+\Big(\dashint_{Q_j} V^{n/2}\Big)^{2/n}<\mu.
\end{equation}

\begin{corollary}[Bounds on the eigenvalue counting function; {\cite[Main Theorem and Remark 4.2]{shenb}}]\label{cor.bounds} Let $n\geq3$, $\Omega=\bb R^n$, and suppose that the operator $L_{\a,0}$ on $\bb R^n$ verifies (\ref{eq.bl1}) and (\ref{eq.ass}). Let $V\geq0$ and assume that $|\B|+V\in RH_{n/2}$. Then there exist $c_1, c_2$, depending only on $n$ and the $RH_{n/2}$ characteristic of $|\B|+V$, so that for every $\mu>0$,
	\begin{equation}\label{eq.eigenbounds}
		\tilde N(c_1\mu)\leq N(\mu,L_{\a,V})\leq\tilde N(c_2\mu).
	\end{equation} 
\end{corollary}

As in the previous results, the main issue to obtain the above corollary without any assumption on $\nabla\B$ was an uncertainty principle (\ref{eq.upm}). Note that only the upper bound is new here, as the lower bound is already proved in \cite{shenb} in a generality that includes our setting. In Remark 4.2 of \cite{shenb}, Z. Shen shows how one can use an uncertainty principle to restore the argument for the main theorem. Then it is an exercise to check that his Lemma 4.1 still works (with proper modifications) in this setting; to check that Lemma 4.1 still works, we mention that it is useful to recall the well-known fact that if $w\in RH_{n/2}$, then $w\in RH_{\frac n2+\ep}$ for some $\ep>0$ determined by the $RH_{n/2}$ characteristic of $w$. Note that we do not need to assume (\ref{eq.bbound}) for this result. 

Since Corollary \ref{cor.bounds}  extends to $n\geq3$ the bounds on the eigenvalue counting function of Shen \cite{shenb} under a ``sign'' assumption on $\B$ (that is, the strong favored directionality assumption (\ref{eq.ass})), while removing the condition (\ref{eq.gradb}) or any other condition on $\nabla\B$, we see that Corollary \ref{cor.bounds} gives a solution  to Problem \ref{pro.s2}.  Example \ref{ex.1} is an explicit magnetic potential which verifies the hypotheses of the corollary.  We also note that Corollary \ref{cor.bounds} is very much related to recent \emph{landscape laws} for non-magnetic Schr\"odinger operators  \cite{dfm4, afmwz}, which give bounds on the integrated density of states in terms of a landscape function.

Next, we show an analogue of \cite[Theorem 0.11]{sheno}.  If $V_-$ is non-trivial, then $L_{\a,V}$ may have negative eigenvalues. The following result gives an upper bound on the number of such negative eigenvalues.

\begin{corollary}[An upper bound on the number of negative eigenvalues of $L_{\a,V}$; {\cite[Theorem 0.11]{sheno}}]\label{cor.neg} Let $n\geq3$, $\Omega=\bb R^n$, and suppose that the operator $L_{\a,0}$ on $\bb R^n$ verifies (\ref{eq.bl1}) and (\ref{eq.ass}) (with $V\equiv0$). Moreover,  let $V\in L^p_{\loc}(\bb R^n)$ for some $p>1$,  and assume that $|\B|\in RH_{n/2}$. Then there exist $C$, which depends only on $n$, and $c,\alpha>0$, which depends only on $n$, $p$,  and the $RH_{n/2}$ characteristic of $|\B|$, such that, for $\mu\leq0$, 
\[
N(\mu,L_{\a,V})\leq CN_0,
\]
where $N_0$ is the number of minimal disjoint dyadic cubes $Q$ which satisfy
\[
\ell(Q)^2\Big(\dashint_Q|V|^p\Big)^{\frac1p}\geq c,\qquad\ell(Q)<\frac1{\sqrt{\mu}},
\]
and
\begin{equation}\label{eq.shencond}
\ell(Q)<\inf_{x\in Q}\frac{\alpha}{m(x,|\B|)}.
\end{equation}
	
\end{corollary}

The condition (\ref{eq.gradb}) in \cite{sheno} was most saliently used to obtain the uncertainty principle (\ref{eq.upm}) and   the decay of the fundamental solution. In our situation, we have already seen that we recover both the uncertainty principle and the exponential decay of the fundamental solution Corollary \ref{cor.fund}. Note that in \cite{sheno}, the number $N_0$ in Theorem 0.11 is defined using the condition
\begin{equation}\label{eq.shencond2}
\ell(Q)^2\Big(\dashint_Q|\B|^2\Big)^{\frac12}\leq1.
\end{equation}
instead of (\ref{eq.shencond}). In fact (see \cite[Remark 0.18]{sheno}), the definition of $N_0$ with (\ref{eq.shencond}) is enough to prove the theorem, and if $\hat N_0$ is defined in the same way as $N_0$ but with (\ref{eq.shencond2}) in place of (\ref{eq.shencond}), then $N_0\lesssim\hat N_0$ whenever we have that
\begin{equation}\label{eq.bbound}
	|\B(x)|\leq C_2m(x,|\B|)^2.
\end{equation}
The inequality (\ref{eq.bbound}) is true if (\ref{eq.gradb}) holds and $|\B|\in RH_{n/2}$, but the novelty of our corollary is that we  replace (\ref{eq.gradb}) with the favored directionality assumption (\ref{eq.ass}). If we further assumed that (\ref{eq.bbound}) is true (which   does not place any assumptions on $\nabla\B$ either), then we could state Corollary \ref{cor.neg} with $\hat N_0$ instead of $N_0$, as it is stated in \cite{sheno}.

Actually,  (\ref{eq.bbound}) is used within the proof of \cite[Theorem 0.11]{sheno} non-trivially, to bound the contribution of a carefully crafted gauge transform ${\bf h}^j$ on a cube $Q_j$ (see \cite[p. 659]{sheno}). However, we may instead obtain the desired bound by following the method of \cite[Lemma 4.1]{shenb}, and using (\ref{eq.shencond}) and the definition of $m(\cdot,|\B|)$, provided that $\alpha$ is small (depending on allowable constants). 

Note we have obtained an appropriate analogue of \cite[Theorem 0.11]{sheno} with no assumptions on $\nabla\B$, and Example \ref{ex.1} shows that Corollary \ref{cor.neg} has content. This resolves Conjecture \ref{conj.s1}.

\subsection{An uncertainty principle  for $L_{\a,V}$ with a good transition region}\label{sec.t} 

Lastly, we remark about the applicability of the landscape function (\ref{eq.fmu}) when $V\geq0$ but the magnetic field does not verify the favored directionality assumption (\ref{eq.signass}) everywhere on $\bb R^n$. As mentioned in the introduction, a scale-invariant regularity condition on the direction of the magnetic field has been used to obtain uncertainty principles, and the condition (\ref{eq.gradb}) of Z. Shen essentially fulfills this role. Let us first see that if (\ref{eq.gradb}) is verified over the whole domain, then we may still construct a landscape function. 

\begin{example}[A landscape function for $L_{\a,V}$ under Shen's assumptions] Suppose  that $\B$ is a magnetic field which is a non-degenerate Shen potential and that verifies the assumption (\ref{eq.gradb}), so that its direction is slowly varying, in particular. In this case, instead of considering the landscape of the operator (\ref{eq.electric}), take $u$ to be instead the landscape function of the operator
	\begin{equation}\label{eq.withb}
		-\Delta +|\B|+V,\qquad\text{so that }-\Delta u+|{\B}|u+Vu=1.
	\end{equation}
	By Theorem \ref{thm.compare}, $\frac1{u(x)}\approx m^2(x,|{\B}|+V)$, and in this case we also have the uncertainty principle \cite[Theorem 2.7]{shene}
	\begin{equation}\label{eq.up3}
		\int_{\bb R^n}\frac1u|f|^2\leq C\Big[\int_{\bb R^n}|D_{\a}f|^2+V|f|^2\Big],\qquad\text{for any }f\in C_c^{\infty}(\bb R^n),
	\end{equation}
	where $C$ depends on $n$, $C_0$, $C_1$, and $\delta$.\hfill{$\square$}
\end{example}

Now let us describe a situation where we may obtain a landscape function and an uncertainty principle over $\bb R^n$ if there is a ``good transition'' region between two subsets of $\bb R^n$ where our stronger directionality assumption (\ref{eq.ass}) is satisfied. Our assumptions are by no means optimal, but the point  is to illustrate the fact that the favored directionality assumptions can be considered alongside with regularity assumptions to establish uncertainty principles.

\begin{theorem}[A landscape for $L_{\a,V}$ with a transition region]\label{thm.transition}  Let $n\geq3$, $\Omega=\bb R^n$, $\B\in RH_{n/2}$ such that $|\B|\geq1$, and take $V\equiv0$. Assume  that $\bb R^n=\Omega_1\cup T\cup\Omega_2$, where $\Omega_1,\Omega_2,T$ are open sets,  
\[
\partial\Omega_i\subset T,\qquad d_i:=\dist(\Omega_i\backslash T,\partial\Omega_i)>0,\qquad i=1,2,
\]
and such that (\ref{eq.ass}) is verified over $\Omega_1$ and over $\Omega_2$ separately\footnote{There may be  different admissible selections $\m S_1$, $\m S_2$ which verify (\ref{eq.magb2}) over $\Omega_1$, $\Omega_2$ respectively, but   the coordinate system is fixed.}. Moreover, assume that there exist $K,M\geq1$ such that  $M\geq\sqrt{\mu}$, where $\mu=\mu(K)$ is fixed in (\ref{eq.mu}), and
\begin{equation}\label{eq.nablab}
|\nabla\B(x)|\leq Km^2(x,|\B|+1),\quad\text{for each }x\in T+M,
\end{equation}
where $T+M:=\{x\in\bb R^n:\dist(x,T)<M\}$. Then
\[
\int_{\bb R^n}\frac1u|f|^2\leq C\int_{\bb R^n}|D_\a f|^2,\qquad\text{for any }f\in C_c^{\infty}(\bb R^n),
\] 
where $u$  is the landscape function for the operator (\ref{eq.withb}) over $\bb R^n$.
\end{theorem}

$T$ should be thought of as a ``transition region'' from $\Omega_1$ to $\Omega_2$ around which we assume that $\nabla\B$ is uniformly well-behaved.
	
\noindent\emph{Proof.}	Fix $f\in C_c^{\infty}(\bb R^n)$. Let $\varphi\in C_c^{\infty}(\Omega_1)$ with $\varphi\equiv1$ on $\Omega_1\backslash T$, $0\leq\varphi\leq1$ on $\Omega_1$, and $|\nabla\varphi|\lesssim d_1^{-1}$. Since $D_{\a}(f\varphi)=\varphi_1D_{\a}f+f\nabla\varphi$, we may obtain  the estimate
	\begin{equation}\label{eq.upest}
		\int_{\Omega_1\backslash T}\frac1u|f|^2\leq C_n\Big[\int_{\Omega_1}|D_{\a}f|^2+\int_{T}d_1^{-1} m^2(\cdot,|\B|)|f|^2\Big],
	\end{equation}
	similar as to how we obtained (\ref{eq.up}), and where we have used that $|\B|\geq1$. Here, $C_n$ is a constant depending only on dimension. We can obtain an analogous estimate with $\Omega_1$, $d_1$ replaced by $\Omega_2$, $d_2$. Now, in the same manner as to how \cite[estimate (2.13)]{shene} is achieved, we may see that
	\begin{equation}\label{eq.t}
		\int_Tm^2(\cdot,|\B|)|f|^2\leq C_2\Big[\int_{T+1}|D_{\a}f|^2+\int_{T+1}\frac{|\B|^2}{m^2(\cdot,|\B|)}|f|^2\Big],
	\end{equation}
	where $C_2$ depends only on $n$ and the $RH_{n/2}$ characteristic of $|\B|$. Next we fix $\phi\in C_c^{\infty}(T+M)$ with $0\leq\phi\leq1$, $\phi\equiv1$ on $T+1$, and $|\nabla\phi|\lesssim1/M$, and thus it may be proved that
	\begin{equation}\label{eq.comm}
		\int_{T+1}\frac{|\B|^2}{m^2(\cdot,|\B|)}|f|^2\leq C_3\Big[\int_{\bb R^n}|D_\a(f\phi)|^2+\int_{\bb R^n}|D_\a(f\phi)||m(\cdot,|\B|)(f\phi)|\Big],
	\end{equation}
	in the same way as \cite[estimate (2.12)]{shene} is shown. Here, $C_3$ depends on $K$ from (\ref{eq.nablab}), $n$, and the $RH_{n/2}$ characteristic of $|\B|$. Putting (\ref{eq.t}) together with (\ref{eq.comm}) and using the Cauchy inequality with a fixed $\ep>0$, we obtain that
	\begin{equation}\label{eq.t2}
		\int_Tm^2(\cdot,|\B|)|f|^2\leq\frac C{\ep}\int_{\bb R^n}|D_\a f|^2+\frac{C_4}{\ep M^2}\int_{\bb R^n}m^2(\cdot,|\B|)|f|^2+\ep\int_{\bb R^n}m^2(\cdot,|\B|)|f|^2.
	\end{equation}
	Bringing the estimates (\ref{eq.upest}), (\ref{eq.t2}) together, and using Theorem \ref{thm.compare}, by choosing $\ep$ appropriately  we may prove that
	\[
	\int_{\bb R^n}m^2(\cdot,|\B|)|f|^2\leq C\int_{\bb R^n}|D_\a f|^2+\frac{C_5}{M^2}\int_{\bb R^n}m^2(\cdot,|\B|)|f|^2+\frac14\int_{\bb R^n}m^2(\cdot,|\B|)|f|^2,
	\]
	where $C, C_5$ depend on $K$, $n$, $d_1$, $d_2$, and the $RH_{n/2}$ characteristic of $|\B|$. If we take
	\begin{equation}\label{eq.mu}
	\mu=4C_5,
	\end{equation}
	 then we finally obtain
	\[
	\int_{\bb R^n}m^2(\cdot,|\B|)|f|^2\leq C\int_{\bb R^n}|D_\a f|^2,
	\]
	which yields the desired estimate by Theorem \ref{thm.compare}.\hfill{$\square$}

\appendix

\section{Existence of Green's function as an integral kernel}\label{sec.green}

In this appendix, we prove that, in the setting of Theorem \ref{thm.exist}, there always exists a Green's function (that is, an integral kernel for the solution operator $L^{-1}$) in the sense of Definition \ref{def.green}. The existence of  Green's function has been considered by several authors   under certain assumptions on the coefficients \cite{lsw, gw, kn85, hk, dhm, simonsemigroups, christ91, shenf, ks1, sak2}, but the purpose of our construction is to give a unified existence proof to the case where $A$ is not necessarily symmetric, $0\leq V\in L^1_{\loc}$ and $n\geq1$.  On the other hand, note that, for our purposes, we are satisfied with the existence of the Green's function, and we do not pursue stronger properties which have been proved for several cases in the aforementioned references, such as continuity of the Green's function or regularity of its gradient.

Our method is classical, relying on  the Dunford-Pettis Theorem, but our very general assumptions on $A$ and $V$ require some technical care, as the solutions to $-\dv A\nabla+V$ may not be continuous; however, the non-negativity assumption on $V$ and the weak maximum principle will  save the day.

\begin{proposition}[Existence of Green's function]\label{prop.green}  Retain the setting of Theorem \ref{thm.exist}. Then the Green's function exists in the sense of Definition \ref{def.green}, and is unique. Moreover, for each bounded set $M\subset\Omega$ and each $p\in[1,\frac{n}{n-2})$ if $n\geq3$ and $p\in[1,\infty)$ if $n=1,2$, we have that
\begin{equation}\label{eq.greenlp}
\sup_{x\in M}\Vert G(x,\cdot)\Vert_{L^p(M)}<+\infty.
\end{equation} 
\end{proposition}  

\noindent\emph{Proof.} Fix $x_0\in\Omega$ and $R_0\in\bb N$ sufficiently large so that $\int_{M_{R_0}}V>0$, $M_{R_0}:=B(x_0,R_0)\cap\Omega$, fix $R>R_0$, and  observe that $M_{R}$ has Lipschitz boundary. We will show that for every $q>\max\{1,n/2\}$ and every $f\in L^q(M_R)$,
\begin{equation}\label{eq.est}
\Vert L^{-1}f\Vert_{L^{\infty}(M_R)}\lesssim\Vert f\Vert_{L^q(M_R)}.
\end{equation}
 
Let $V_N:=\min\{V,N\}$ for each $N\in\bb N$. Let $0\leq f\in L^q(M_R)$ for some $q>\max\{1,n/2\}$, and for each $N\in\bb N$, write $\psi_N=L^{-1}_{V_N}f$. By Lemma \ref{lm.approx}, we have that $\psi_N\searrow\psi:=L_V^{-1}f$, and therefore by the global boundedness of weak solutions on bounded domains \cite[Theorem 8.25]{gtpde}\footnote{If $n=1$, the same desired result will follow from Morrey's inequality.}, we see that
\begin{equation}\nonumber
	\sup_{M_R}\psi\leq\sup_{M_R}\psi_1\lesssim\Vert\psi_1\Vert_{L^2(M_{2R})}+\Vert f\Vert_{L^q(M_R)}\lesssim\Vert\psi_1\Vert_{\f l}+\Vert f\Vert_{L^q(M_R)}\lesssim\Vert f\Vert_{L^q(M_R)},
\end{equation}
where we have used (\ref{eq.coercive}) and the Sobolev embedding, and the implicit constant depends only on $n$, $\lambda$, $q$, $R$, and $M_R$. Thus we have shown (\ref{eq.est}) for  any $0\leq f\in L^q(M_R)$. Then we can obtain the general case  $f\in L^q(M_R)$ by the weak maximum principle and using (\ref{eq.est}) for $|f|$ (consider $L^{-1}(|f|-f)$, and as before, approximations via bounded potentials); we omit further details. 

Thus the operator $L^{-1}$ maps $L^q(M_R)$ into $L^{\infty}(M_R)$\footnote{Or rather, a restriction of $L^{-1}$ to $M_R$ does.}. By the Dunford-Pettis Theorem (see \cite[Theorem 46.1]{treves} or \cite[Theorem 2]{bgp07}), it follows that $L^{-1}$ has a unique integral kernel on $M_R$, and so   there exists a measurable function $G_R:M_R\times M_R\ra\bb R$ such that
\begin{equation}\label{eq.green21}
	(L^{-1}f)(x)=\int_{M_R}G_R(x,y)\,f(y)\,dy,\qquad\text{for a.e. }x\in M_R,
\end{equation}
and if $q'$ is the H\"older conjugate of $q$, then
\[
\sup_{x\in M_R}\Vert G_R(x,\cdot)\Vert_{L^{q'}(M_R)}<+\infty.
\]
Observe that if $R'>R$, then for every $x\in M_R$, $G_{R'}(x,\cdot)$ must coincide pointwise a.e. on $M_R$ with $G_R(x,\cdot)$. From these facts, we deduce that there exists a measurable function $G:\Omega\times\Omega\ra\bb R$ such that for any $q>n/2$ and $f\in L^q_c(\Omega)$, the identity (\ref{eq.green}) holds, and furthermore, for each $x\in\Omega$ we have that $G(x,\cdot)\in L^p_{\loc}(\Omega)$ for any $p\in[1,\frac n{n-2})$. The non-negativity follows immediately from the maximum principle.\hfill{$\square$}

\bibliography{../References/refs-nobadurls} 
\bibliographystyle{alphaurl-5} 

\end{document}